\documentclass[a4paper, 11pt]{amsart}

\usepackage{amsmath, amssymb, amsthm}
\usepackage[english]{babel}
\usepackage[utf8]{inputenc}
\usepackage{enumerate}
\usepackage{xcolor}
\def\bx{{\bf x}}
\def\bX{{\bf X}}

\def\by{{\bf y}}
\def\bm{{\bf m}}

\def\bv{{\bf v}}

\newcommand{\ee}{\mathbf{e}}

\def\e{\epsilon}

\def \be {  \varpi}

\newcommand{\fer}[1]{(\ref{#1})}

\def\MI{ {\mathcal {I}}}

\newcommand{\R}{\mathbb R}
\def\be#1\ee{\begin{equation}#1\end{equation}}

\newcommand{\bq}{\begin{equation}}
\newcommand{\eq}{\end{equation}}

\newtheorem{thm}{Theorem}

\newtheorem{prop}[thm]{Proposition}
\theoremstyle{remark}
\newtheorem{rem}[thm]{Remark}
\theoremstyle{definition}

 \numberwithin{equation}{section}

\usepackage{verbatim}

\newenvironment{equations}{\equation\aligned}{\endaligned\endequation}

\begin{document}

\title[Fokker--Planck equations and inequalities for isotropic densities]{Fokker--Planck equations and n--dimensional Poincaré inequalities for isotropic densities}
%
\author{GIULIA FURIOLI}
\thanks{DIGIP, University of Bergamo, viale Marconi 5, 24044 Dalmine, Italy\\
giulia.furioli@unibg.it}
\author{ADA PULVIRENTI}
\thanks{Department of Mathematics, University of Pavia, 
via Ferrata 1,
Pavia, 27100 Italy\,\,\,\\
ada.pulvirenti@unipv.it}
\author{ELIDE TERRANEO}
\thanks{Department of Mathematics,
University of Milan, via Saldini 50, 20133 Milano, Italy \\
elide.terraneo@unimi.it}
\author{GIUSEPPE TOSCANI}
\thanks{Department of Mathematics,, University of Pavia,
via Ferrata 1,
Pavia, 27100 Italy,\\
giuseppe.toscani@unipv.it}

\maketitle

%
%
\begin{center}\small
\parbox{0.85\textwidth}{

\textbf{Abstract.} 
We consider new connections between the problem of trend to equilibrium for the n-dimensional Fokker--Planck equation of statistical physics, and weighted Poincar\'e inequality. To this aim we consider a class of n-dimensional Fokker--Planck equations with variable isotropic coefficient of diffusion and drift, inspired by the analogous one-dimensional Fokker--Planck equation appearing when studying the evolution of wealth distribution.}
\end{center}

\medskip


\medskip


\section{Introduction}

This paper deals with n-dimensional functional inequalities of the type of  Poincar\'e, with weight, for \emph{isotropic} probability densities. Unlike previous studies \cite{BL}, the derivation of these inequalities is strongly related to their connection with the problem of convergence to equilibrium for a class of Fokker--Planck type equations characterized by a variable coefficient of diffusion.  

Let $\bX =(X_1,\dots,X_n)$ be a random vector with $n \ge1$ components, distributed with probability density $f(|\bx|)$, where  $\bx= (x_1,\dots, x_n)\in \MI \subseteq \R^n$. The random vector $\bX$ is said to satisfy a weighted Chernoff--Poincar\'e--type inequality with weight function $w(\bX)$ (where $w$ is a fixed nonnegative,  Borel measurable function), if for any bounded smooth scalar function $\phi$ on $\MI$ 
 \be\label{Poi}
 Var\left[\phi(\bX)\right] \le E\left\{w(\bX)|\nabla \phi (\bX)|^2\right\}.
 \ee
 As usual, for a given scalar function $\phi$
 \[
 E(\phi(\bX)) = \int_\MI \phi(\bx)f(|\bx|)\, d\bx,
 \]
 denotes its expectation value, and
 \[
 Var\left[\phi(\bX)\right] = 
 \int_{\MI} \left( \phi(\bx)- \int_{\MI} \phi(\by)\, f(|\by|) \, d\by \right)^2  f(|\bx|) \, d\bx
 \]
denotes its variance. \\
This inequality is understood in the following sense: if the right-hand side is finite, then the inequality holds true.
For brevity, in the following we will denote the inequality as a Poincaré inequality, the Chernoff denomination being mostly used in a probabilistic framework.

The relationships between the Fokker--Planck equation and differential inequalities like \fer{Poi} has a long history.\\
In divergence form, the classical Fokker--Planck equation 
reads
\be
\label{FP-n}
 \frac{\partial f(\bx, t)}{\partial t} = \nabla\cdot\left(\nabla f(\bx,t) +  f(\bx,t)\nabla\Gamma(\bx)  \right).
\ee
 The unknown $f = f(\bx,t)$  stands for the density of a system of particles at time $t\ge 0$, and without loss of generality can be assumed to be a probability distribution on $\R^n$, $n\ge 1$, since equation \fer{FP-n} conserves both non-negativity and the mass of the solution. In statistical physics, the Fokker--Planck equation models a set of particles subject to a constant diffusion and a drift characterized by a potential $\Gamma(\bx)$, and the balance between these two processes determines  the shape of the equilibrium density
\be\label{eq-n}
f_\infty(\bx) = C\exp\{ -\Gamma(\bx)\},
\ee
which is the (unique) probability density satisfying the differential equation
\be\label{flux}
\nabla f(\bx) +  f(\bx)\nabla \Gamma(\bx) = 0.
\ee
Equation \fer{flux}, identifying the equilibrium density \fer{eq-n} easily follows by applying the divergence theorem to the Fokker--Planck equation \fer{FP-n} on a volume $V \subset  \R^n$ bounded by a closed surface $S$. In fact 
\begin{equations}\label{dive}
&\frac d{dt} \int_V f(\bx,t)\, d\bx = \int_V  \frac{\partial f(\bx, t)}{\partial t}\, d\bx = \\
& \int_V \nabla\cdot\left(\nabla f(\bx,t) +  f(\bx,t)\nabla \Gamma(\bx)\right) \, d\bx = \int_S \left(\nabla f(\bx,t) +  f(\bx,t)\nabla \Gamma(\bx)\right) \cdot {\bf n}\ dS. 
\end{equations}
Thus, whenever $f$ does not depend on time, identity \fer{dive} is satisfied independently of the choice of the volume if and only if $f$ satisfies \fer{flux}. \\
In kinetic theory of rarefied gases, the space variable is the space of velocities $\bv$,  and the potential $\Gamma(\bv)$ is  the kinetic energy $|\bv|^2/2$ \cite{To99}. In this situation the equilibrium density of equation \fer{FP-n}, which reads now
\[
 \frac{\partial f(\bv, t)}{\partial t} = \nabla\cdot\left(\nabla f(\bv,t) + \bv f(\bv,t) \right),
\]
 is the Maxwellian function
\be\label{Max}
M(\bv) = \frac 1{(2\pi)^{n/2} }\exp \{ -|\bv|^2/2\}.
\ee
As it happens in classical kinetic theory for the spatially homogeneous Boltzmann equation, the solution to the Cauchy problem of the Fokker--Planck equation is expected to converge to its equilibrium state, and one would like to estimate the rate of convergence in terms of the initial value.

In \cite{To99}, it has been shown that, in the case  of the squared potential $\Gamma(\bx)= |\bx|^2/2$, the solution converges towards the Maxwellian function \fer{Max} at exponential rate, by resorting to the study of the evolution of the relative entropy functional 
\be\label{ent}
H(f|f_\infty)(t) = \int_{\R^n} f(\bx,t)\log \frac{f(\bx,t)}{f_\infty(\bx)} \, d\bx
\ee
and making use of the logarithmic Sobolev inequality by Gross \cite{Gross}. Moreover, in \cite{To99} it has been pointed out that the logarithmic Sobolev inequality is deeply connected with the entropy dissipation of the Fokker--Planck equation, so that the equation can be used to obtain a physically based proof  of the inequality \cite{To97}.  This way of using evolution equations to get sharp inequalities was already present  in the Bakry--Emery work \cite{BE}, which takes a rather abstract point of view, based on the so-called $\Gamma_2$, or \emph{carr\'e du champ it\'er\'e}.

The importance of the Fokker--Planck type equation \fer{FP-n}  in the derivation of functional inequalities is nowadays well known. The interplay between this kinetic equation  and functional analysis has been explored  about one quarter of century ago in \cite{AMTU,OV}, and excellently reviewed in \cite{MV}, where new connections were brought to light.
Following this line of research, a large amount of the mathematical studies of the last years have been devoted to study the relationships between the rate of convergence of the solution to the Fokker--Planck equation \fer{FP-n}  to its equilibrium, as given by \fer{eq-n}, and the shape of the confinement potential $\Gamma(\bx)$. \\
In particular, starting from the work \cite{TV}, research was directed toward a precise estimate of  the rate of convergence to equilibrium of the solution to the Fokker--Planck equation in presence of a so called {\it subcritical confinement} (cf. the recent paper \cite{KMN} and the references therein).  

As a matter of fact, equilibria of type \fer{eq-n}, which are relative to a Fokker--Planck equation of type \fer{FP-n}, characterized by a constant diffusion and a general potential $\Gamma(\bx)$, can  sometimes appear as equilibria of  Fokker--Planck type equations with quadratic potential and variable diffusion coefficient.  Fokker--Planck type equations with this structure have been recently considered in the modeling  of social and economic phenomena,  a challenging research activity in the communities of both physicists and applied mathematicians, who classified the fields of research with the names of socio-physics and, respectively, econophysics \cite{FPTT, NPT, PT13}. 

One of the typical features of these phenomena is related to the tails of the underlying steady distributions, which are often characterized by polynomial decay at infinity \cite{New}.\\
A prominent example appears in the study of  the distribution of wealth among trading agents \cite{BM, CoPaTo05}. In this economic context, the evolution of the density  $f(x,t)$ of a system of agents with personal wealth $x\in \R, x \ge0$ at time $t \ge 0$ is shown to be solution to  the  Fokker--Planck type equation 
  \be\label{FPw}
 \frac{\partial f(x, t)}{\partial t} = \frac \sigma{2}\frac{\partial^2 }
 {\partial x^2}\left( x^2 f(x,t)\right) + \lambda \frac{\partial }{\partial x}\left(
 (x-1) f(x,t)\right),
 \ee
characterized by a variable coefficient of diffusion and linear drift (or quadratic potential). In \fer{FPw},  $\sigma$ and $\lambda$ denote  positive constants related to essential properties of the trade rules of the agents.  
By fixing the mass density equal to unity, the unique steady state  of equation \fer{FPw} is the inverse Gamma density 
 \be\label{equi2}
f_\infty(x) =\frac{\mu^{1+\mu}}{\Gamma(1+\mu)}\frac{\exp\left(-\frac{\mu}{x}\right)}{x^{2+\mu}}, \quad  \mu = \frac{2\lambda}\sigma
 \ee
 where $\Gamma$ is the Euler Gamma function.
 This stationary distribution, in agreement with the analysis of the Italian economist Vilfredo Pareto \cite{Par}, exhibits a power-law tail for large values of the wealth  variable.
 
 Equation \fer{FPw} is  a one-dimensional example of the class of Fokker--Planck  type equations 
 \be\label{FP-k}
\frac{\partial f(\bx,t)}{\partial t} = \nabla\cdot\left[\nabla (K(\bx) f(\bx,t)) + (\bx-\bm) f(\bx,t) \right], \quad \bx, \bm \in \MI\subseteq \R^n,
\ee
characterized by a variable diffusion coefficient   $K = K (\bx)$ and a quadratic potential.\\
In a bounded domain $\MI\subseteq \R^n$, equation \fer{FP-k} is usually complemented  with the so-called no-flux  boundary conditions
\be\label{bou}
 \nabla (K(\bx)f(\bx,t)) + (\bx-\bm)f(\bx,t) \big|_{\bx \in \partial \MI}= 0, \quad t>0,
\ee
that ensure mass conservation.

In one dimension, functional inequalities related to Fokker--Planck equations of type  \fer{FP-k} have been obtained  recently in \cite{FPTT21}. The analysis in \cite{FPTT21} leads in many cases to inequalities with sharp constants. While the method of proof, based on physical arguments, shares a new light on the meaning of these inequalities, the quality of mathematical results is comparable to other approaches \cite{BCG,BL,BJ,BJM1,BJM2,CGGR,Goz}. 

Let us come now to the content of the paper. In Section \ref{sec:FP} we detail some properties of the Fokker--Planck type equations \fer{FP-k} which clarify both the presence of the weight, and the importance of Poincar\'e--type inequalities with weight in the study of the rate of convergence to equilibrium. In Section \ref{sec:Peso} we introduce various typical examples of isotropic n-dimentional densities and identify the possible weight we expect in a Poincaré--type inequality for these densities. In Section \ref{sec:poincare} we  prove first an extension to the n-dimensional setting of the classical Poincaré inequality obtained in \cite{FPTT}, valid for the class of n-dimensional probability densities expressed as a product of one-dimensional ones (Theorem \ref{Ch_Ind}). In particular,  Subsection \ref{sec:chernoff} 
 deals with an application of the result for product densities to general isotropic probability densities (Theorem \ref{Cher}).
  In Subsection \ref{sec:examples} we apply the previous results  to the densities introduced in Section \ref{sec:Peso}. Unfortunately,  the weight we obtain appears  to be sharp for the Gaussian and generalized Gaussian densities but this is not the case for exponential--type and  Barenblatt densities. Therefore, in Section \ref{sec:refined} we prove a refined
Poincaré--type result (Proposition \ref{C}),  with the expected  weight,  but which only applies to bounded smooth functions  supported outside  a ball. Both  in Section  \ref{sec:poincare} and \ref{sec:refined},  we 
   take   advantage of the representation of isotropic  densities as steady solutions to Fokker--Planck type equations as \fer{FP-k}.
 We conclude our paper with Section \ref{sec:hybrid}. There we  prove  a hybrid weighted inequality (Theorem \ref{H}) in which we get the expected weight but we pay the prize of an extra term, and so this is not a Poincaré--type inequality.

\section{Basic properties of Fokker--Planck type equations}\label{sec:FP}

Before entering into the details of the various proofs, let us briefly explain the strategy behind our approach. In the rest of this section we perform formal calculations and do not deal with their rigorous justifications. We address to the work \cite{LL1} for technical details about  regularity issues. \\
To start with, we recall that both the Fokker--Planck equations \fer{FP-n} and \fer{FP-k} can be rewritten in a form which is more suitable to study the decay of their solutions towards the equilibrium density in relative entropy. For the sake of simplicity, we will refer our analysis to equation \fer{FP-k}, in which the coefficient of diffusion $K$ satisfies $K(\bx) >0$, $\bx \in \MI\subseteq \R^n$.

As already discussed in the introduction, the stationary solution $f_\infty$ of equation \fer{FP-k} is found by solving on $\MI\subseteq \R^n$ the differential equation
 \be\label{stazio}
\nabla\left(K(\bx)
 f_\infty \right) + (\bx -\bm)\,f_\infty = K(\bx) f_\infty \left[ \nabla \log \left(K(\bx) f_\infty \right) + \frac{\bx -\bm}{K(\bx)}\right] = 0, 
 \ee 
 that leads to the identity
 \be\label{id1}
 \nabla \log \left(K(\bx) f_\infty \right) = -  \frac{\bx -\bm}{K(\bx)}.
 \ee
Making use of identity \fer{id1}, one easily shows that the right-hand side of equation \fer{FP-k} can be written in different ways \cite{FPTT}. Indeed
 \[
\nabla \left(K(\bx)
 f \right) + (\bx -\bm)\,f = K(\bx) f \left[ \nabla \log(K(\bx) f) + \frac{\bx-\bm}{K(\bx)}\right]=
 \]
 \[
  K(\bx) f \left[ \nabla \log(K(\bx) f)- \nabla \log(K(\bx) f_\infty) \right] =   K(\bx) f  \, \nabla \log\frac f{f_\infty}= K(\bx) f_\infty  \nabla \frac f{f_\infty}.
 \]
Hence, the Fokker--Planck equation \fer{FP-k} can be rewritten in the equivalent form
 \be\label{FPalt}
  \frac{\partial f}{\partial t} = \nabla\cdot \left[ K(\bx) f \,\nabla \log\frac f{f_\infty}\right],
 \ee
which enlightens the role of the logarithm of the quotient $f/f_\infty$, and
 \be\label{FPal2}
  \frac{\partial f}{\partial t} = \nabla\cdot \left[ K(\bx) f_\infty \nabla \frac f{f_\infty}\right].
 \ee 
In particular, owing to \fer{stazio}, the form \fer{FPal2} {allows us} to obtain the evolution equation for the quotient $F= f/f_\infty$. 
We have
\begin{equations}\nonumber
 \nabla\cdot \left[ K(\bx) f_\infty \nabla F \right] &= K(\bx) f_\infty \Delta F + \nabla\left(K(\bx) f_\infty \right)\cdot \nabla F = \\
 &K(\bx) f_\infty \Delta F -f_\infty (\bx-\bm) \cdot \nabla F.
\end{equations}
Therefore
 \[
  \frac{\partial f}{\partial t} = f_\infty  \frac{\partial F}{\partial t} = K(\bx) f_\infty \Delta F -f_\infty (\bx-\bm) \cdot \nabla F,
 \]
which shows that $F$ satisfies the (adjoint) equation
 \be\label{quo}
\frac{\partial F}{\partial t} =  K(\bx) \Delta F  -(\bx-\bm) \cdot\nabla F.
 \ee
In analogous way, the boundary conditions of the two equivalent forms   \fer{FPalt} and \fer{FPal2} of the Fokker--Planck equation \fer{FP-k} follow by imposing mass conservation.  If \fer{FPalt} is concerned, the natural boundary condition is given by
\be\label{BCalt} 
 K(\bx) f(\bx,t) \nabla \log\frac {f(\bx,t)}{f_\infty(\bx)} \Big|_{\bx\in \partial\MI } = 0, \quad t >0.
\ee 
Likewise, if mass conservation is imposed on equation \fer{FPal2},  the natural boundary condition reads 
\be\label{BCal2} 
K(\bx) f_\infty(\bx) \nabla \frac {f(\bx,t)}{f_\infty(\bx)}  \Big|_{\bx\in \partial\MI } = 0, \quad t >0.
\ee 
For the sake of simplicity
 in the rest of this section we limit ourselves  to  spherical domains. Consequently $\MI = B_R= \{ \bx : |\bx| \le R\}$,  $0 < R \leq \infty$, and $\partial \MI = \{ \bx : |\bx| = R\}$ with obvious meaning for $R=\infty$.  
In presence of no-flux boundary conditions, and enough regularity of the solution to the Cauchy problem, one can easily extend two one-dimensional results proven in \cite{FPTT} to the n-dimensional  Fokker--Planck equation \fer{quo}, which relate Poincaré inequality to the convergence to the equilibrium for the Fokker--Planck equation  \fer{FP-k}.
They both allow us to obtain convergence in $L^1$ but the price to pay in each case is different: in the first case, we assume that the initial data is very close to the stationary state itself and we obtain exponential convergence, while in the second case we assume less restriction on the initial date at the price of a slower (polynomial) convergence. Here are the results.
\begin{thm}\label{thm1}\cite{FPTT} Let $f_0$  be  the initial data of a Cauchy problem for equation  \fer{FP-k} coupled with  boundary conditions \fer{bou}
and let us assume
\be\label{cond-1}
  \int_\MI \left[\left( \frac{f_0(\bx)}{f_\infty(\bx)}\right)^2 - 1 \right] f_\infty(\bx) \, d\bx <+\infty
\ee
where $f_\infty$ is the stationary state and satisfies the Poincaré inequality \fer{Poi}  with weight $w(|\bx|) = c K(|\bx|)$, with $c$ positive constant. 
Then the solution $f(t)$ converges to $f_\infty$ at an exponential rate and 
\[
\int_\MI \left[\left( \frac{f(\bx,t)}{f_\infty(\bx)}\right)^2 - 1 \right] f_\infty(\bx) \, d\bx \leq e^{-\frac 2c t}   \int_\MI \left[\left( \frac{f_0(\bx)}{f_\infty(\bx)}\right)^2 - 1 \right] f_\infty(\bx) \, d\bx.
\]
\end{thm}
\begin{rem} Since 
\[
\begin{aligned}
\int_\MI |f(\bx, t)-f_\infty(\bx)|\, d\bx  & \leq  \left( \int_\MI  \frac {\left (f(\bx,t)-f_\infty (\bx)\right )^2} {f_\infty(\bx)} \, d\bx\right )^{1/2} { \left(\int_\MI f_\infty(\bx)\, d\bx\right )^{1/2}}\\
&= \left(\int_\MI \left[\left( \frac{f(\bx,t)}{f_\infty(\bx)}\right)^2 - 1 \right] f_\infty(\bx) \, d\bx\right ) ^{1/2}
\end{aligned}
\]
we easily obtain exponential convergence in $L^1$ for such initial data.
\end{rem}
The second result involves the Hellinger distance 
\[
d_H(f( t), f_\infty)= \left[\int_\MI \left(\sqrt {f (\bx,t)} -\sqrt {f_\infty(\bx) }\right )^2\ d \bx\right]^{1/2} 
\]
 between the solution and the stationary state and reads as follows.
 \begin{thm}\label{thm2}\cite{FPTT}
 Let $f(t)$ be the solution  of the Fokker--Planck equation \fer{FP-k}, corresponding to an initial value $f_0$ such that the relative Shannon entropy $H(f_0,f_\infty) = \int_{\MI} f_0(\bx) \log \frac {f_0(\bx)}{f_\infty(\bx)}\, d \bx$ and the Hellinger distance $d_H(f_0, f_\infty)$ are bounded and coupled with  boundary conditions \fer{bou}, where  $f_\infty$ is the stationary state. Let us assume that $f_\infty$  sa\-tis\-fies the Poincaré inequality \fer{Poi}  with weight $w(|\bx|) = c K(|\bx|)$, with $c$ positive constant.  
 Then,
\[
d_H(f(\bx,t), f_\infty(\bx))
  = o\left(\frac 1{\sqrt t} \right), \quad t\to \infty.
\]
 \end{thm}
 
 {\begin{rem} The decay of the Hellinger distance implies decay in $L^1(\MI)$, since
 \[
\begin{aligned}
&\int_\MI |f(\bx,t)-f_\infty(\bx)| \ d\bx = \int_\MI \left|\sqrt {f(\bx,t)}-\sqrt {f_\infty(\bx)}\right | \left(\sqrt{f(\bx,t)} + \sqrt {f_\infty(\bx)}\right) \ d\bx\\
&\leq  \left(\int_\MI \left(\sqrt {f(\bx,t)}-\sqrt {f_\infty(\bx)}\right)^2 \ d\bx\right )^{\frac 12} \left( \int_\MI \left(\sqrt{f(\bx,t)} + \sqrt {f_\infty(\bx)}\right)^2 \ d\bx\right )^{\frac 12}\\
&= d_H(f(t), f_\infty)  \left( \int_\MI \left(f(\bx,t) + f_\infty(\bx) + 2 \sqrt {f(\bx,t) f_\infty(\bx)} \right) \ d\bx\right )^{\frac 12}\\
& = \sqrt 2 d_H(f(t), f_\infty) \left(1+ \int_I  \sqrt {f(\bx,t) f_\infty(\bx)}\ d\bx \right )^{\frac 12}
\end{aligned}
\]
as 
\[
\int_\MI  \sqrt {f(\bx,t) f_\infty(\bx)} \ d\bx 
\leq 
\left( \int_\MI( \sqrt {f(\bx,t)})^2 \ d\bx \right )^{\frac{1}{2}}
\left( \int_\MI(\sqrt {f_{\infty}(\bx)} )^2  \ d\bx \right )^{\frac{1}{2}}\\
\]
which, in turns, implies
\[
\|f(t)-f_\infty\|_{L^1} \leq 2  d_H(f(t), f_\infty).
\]
This leads to the polynomial convergence of the solution to the Fokker--Planck type equations \fer{FP-k}  towards equilibrium in presence of natural conditions on the initial value (finite relative entropy) any time a Poincar\'e inequality for the equilibrium density is shown to hold.
\end{rem}

The two results are consequence of the following Theorem, which is a $n$--dimensional version of the analogous Theorem 3.1 in \cite{FPTT}.

\begin{thm}\cite{FPTT} \label{giusto} 
Let the smooth function $\phi(r)$, $r\in \R_+$ be convex. Then, if $F(\bx,t)$ is the solution to equation \fer{quo} in $\MI=  \{ \bx : |\bx| \le R\} $, $0 < R \leq \infty$, satisfying  the boundary condition \eqref{BCal2} and $c \le F(\bx,t) \le C$ for some positive constants $c<C$, then  the functional
\[
\Theta(F(t)) = \int_\MI f_\infty(\bx) \phi(F(\bx,t))\, d\bx 
\]
is monotonically decreasing in time, and the following equality holds true
 \be\label{dec4}
 \frac d{dt} \Theta (F(t)) = -  I_\Theta (F(t)), 
 \ee
 where $I_\Theta$ denotes the nonnegative quantity
 \be\label{ep1}
 I_\Theta (F(t)) =  \int_\MI K(\bx) f_\infty(\bx) \left|\nabla F(\bx,t) \right|^2 \phi''(F(\bx,t))  \, d\bx
 \ee
 for $f_\infty$ the stationary state of equation \eqref{FP-k}.
\end{thm}
The proof of the theorem can be found in the Appendix and follows the same lines as the analogous one  in \cite{FPTT}.

{\it Proof of Theorem \ref{thm1}}.  Theorem \ref{thm1} follows by choosing in Theorem  	\ref{giusto} the convex function $\phi(r) = (r-1)^2$, which implies $\phi''(r)= 2$. Then
\[
\Theta(F(t)) = \int_\MI f_\infty(\bx)\left[ F(\bx,t)-1 \right]^2\, d\bx, 
\]
and
\[
I_\Theta (F(t)) = 2 \int_\MI K(\bx) f_\infty(\bx) \left|\nabla F(\bx,t) \right|^2  \, d\bx.
\]
Suppose now that the probability density function $f_\infty$ satisfies the Poincar\'e inequality \fer{Poi}  with weight $w(|\bx|) = c K(|\bx|)$, with $c$ positive constant. Then, if at any given time $t>0$, we choose $\phi(\bx)= {f(\bx,t)}/{f_\infty(\bx)}$, provided $\Theta(F(t=0))< +\infty$ we obtain
\begin{equations}\label{poi1}
I_\Theta (F(t)) &\ge  \frac 2c \int_\MI \left[\left( \frac{f(\bx,t)}{f_\infty(\bx)}\right)^2 - \left( \int_\MI f_\infty(\bx)\frac{f(\bx,t)}{f_\infty(\bx)} \, d\bx \right)^2 \right] f_\infty(\bx) \, d\bx= \\
&  \frac 2c \int_\MI \left[\left( \frac{f(\bx,t)}{f_\infty(\bx)}\right)^2 - 1 \right] f_\infty(\bx) \, d\bx = \frac 2c \, \Theta(F(t)).
\end{equations}
Coupling \fer{poi1} with \fer{dec4} we conclude with the exponential convergence of the solution to the Fokker--Planck equation \fer{FP-k} towards equilibrium at exponential rate in a weighted $L^2$-setting.  
\hfill $\square$

 \medskip
 {\it Proof of Theorem \ref{thm2}}. 
Consider now in Theorem  \ref{giusto} the classical case in which the convex functional $\Theta (F(t))$ is the relative entropy,
associated to the convex function $\phi(r) = r\log r$
so that
\[
\Theta (F(t)) = \int_\MI f_\infty(\bx) F(\bx,t)\log F(\bx,t) \, d\bx = \int_\MI f(\bx,t) \log \left( \frac {f(\bx,t)}{f_\infty(\bx)}\right ) d \bx.
\]
In this case 
\begin{equations}\label{ent1}
 I_\Theta (F(t)) &=  \int_\MI K(\bx) f_\infty(\bx) \left|\nabla F(\bx,t) \right|^2 \frac 1{F(\bx,t)}  \, d\bx =  \\
&4 \int_\MI K(\bx) f_\infty(\bx) \left|\nabla F(\bx,t) \frac 1{2\sqrt{F(\bx,t)}}\right|^2   \, d\bx = \\
&4 \int_\MI K(\bx) f_\infty(\bx) \left|\nabla \sqrt{F(\bx,t)}\right|^2   \, d\bx.
 \end{equations}
 If we now apply Poincar\'e  inequality \fer{Poi}  with $\phi(\bx) =\sqrt{{f(\bx,t)}/{f_\infty(\bx)}}$ at fixed time $t>0$, we  obtain
\be\label{formula}
\begin{aligned}
& I_\Theta (F(t)) =4 \int_{\MI}  K (\bx)  f_\infty  (\bx)
\left|\nabla\sqrt{\frac {f (\bx,t)}{f_\infty (\bx)}}\right|^2\ d\bx \geq  \\
&\frac 4c\left[ \int_\MI \frac {f (\bx,t)}{f_\infty (\bx)} f_\infty (\bx) \ d\bx -   \left( \int_\MI\sqrt {\frac {f (\bx,t)}{f_\infty (\bx)}} f_\infty(\bx)\ d\bx\right )^2\right]= \\
&
\frac 4c \left[ 1 -  \left(\int_\MI \sqrt{ {f (\bx,t)}\, {f_\infty(\bx) }} \ d\bx\right )^2\right].
\end{aligned}
\ee
The lower bound in \fer{formula} can be easily related to the Hellinger distance \cite{FPTT}.
\noindent Indeed, for any pair of probability density functions $f(\bx)$ and $g(\bx)$ supported in $\MI$ we have
\be\label{hell-dis}
\begin{aligned}
&\int_\MI \left(\sqrt {f (\bx)} -\sqrt {g(\bx) }\right )^2\ d \bx  = \int_\MI\left( f (\bx) + g(\bx)  - 2\sqrt{f(\bx)\, g(\bx)}\right ) d\bx =\\
& 2\left(1-\int_\MI \sqrt {f(\bx)\, g(\bx)}  \ d\bx \right ) \leq 2\left(1-\left(\int_\MI \sqrt {f(\bx)\, g(\bx)}\ d\bx \right )^2\right ).
\end{aligned}
\ee
The last inequality in \fer{hell-dis} follows from  Cauchy--Schwartz inequality, as $\int_\MI \sqrt {f(\bx)\, g(\bx)} d\bx\leq 1.$
Thus, starting from the time-decay of the relative entropy,  we conclude with the bound
\be\label{con33}
I_\Theta (F(t))  \geq  \frac{2 }{c}\int_\MI \left(\sqrt {f (\bx,t)} -\sqrt {f_\infty(\bx) }\right )^2\ d \bx = \frac{2 }{c}\left(d_H(f(t),f_\infty)\right )^2.
\ee
Therefore, if we assume that the relative entropy is initially bounded, by \fer{dec4} we can conclude with the integrability on time 
of the square of the Hellinger distance
\be\label{final}
\int_0^\infty d_H(f(s),f_\infty)^2 \, ds \le \frac c2 H(f_0, f_\infty) .
\ee
Now, since the function $\phi(r) = (\sqrt r -1)^2$ is convex, Theorem \ref{giusto} applies again, thus giving that the square of the Hellinger distance  between the solution to the Fokker--Planck equation \fer{FP-k} and its equilibrium density is monotonically decresing in time (see \cite{FPTT} for further details). Coupling integrability condition \fer{final} with the monotonicity of $d_H(f(s),f_\infty)^2$ we can conclude with the time decay of the Hellinger distance at a rate almost of $o(1/\sqrt t)$.  
  \hfill $\square$

\section{Fokker--Planck equilibria and weight in Poincar\'e--type inequalities}\label{sec:Peso}

It is interesting to remark that, while the knowledge of the diffusion function $K(\bx)$ allows us to determine in a unique way the equilibrium density  $f_\infty$ of the Fokker--Planck equation \fer{FP-k} by solving the differential equation \fer{id1}, any isotropic density function is  uniquely related to a coefficient of diffusion $K(|\bx|)$ in the following way. Let $f_\infty(|\bx|)$ with $x\in \MI= \{   |\bx| \leq i_+ \}$ be an isotropic probability density function which represents the steady state solution of the Fokker--Planck equation \fer{FP-k}, where $\mathbf {m=0}$, namely
 \be\label{FP-k}
\frac{\partial f(\bx,t)}{\partial t} = \nabla\cdot\left[\nabla (K(\bx) f(\bx,t)) + \bx f(\bx,t) \right], \quad \bx \in \MI
\ee
 and let us suppose  $f_\infty >0$ and such that $\int_0^{i_+} \rho f_\infty(\rho) d\rho <+\infty$.
In this case, it is immediate to conclude that the exact value of the coefficient of diffusion $K(\bx)$ can be expressed by the formula
\be\label{K-ok}
K(\bx) = K(|\bx|)=  \frac{\int_{|\bx|^2}^{(i_+)^2} f_\infty(\sqrt y)\, dy}{2 f_\infty(|\bx|)}.
\ee 
In fact, \fer{K-ok} implies that $f_\infty$ is a steady state for \eqref{FP-k}.  Hence, thanks to the discussion of Section \ref{sec:FP} we can conjecture that $w(|\bx|) = c K(|\bx|)$,  where $K(\cdot)$ is defined by \fer{K-ok}, is the right weight  to get the Poincar\'e inequality \fer{Poi}. \\
Definition \fer{K-ok} establishes an easy-to handle connection between a probability density and its weight in Poincar\'e inequality \fer{Poi}, and an interesting connection between steady states solutions of Fokker--Planck type equations like \fer{FP-k} and Poincar\'e inequalities with weight.\\
Few examples will clarify this relationship. Let $f_\infty(|\bx|)$ be a centered Gaussian density in $\R^n$, i.e.
\be\label{Normal}
f_\infty(|\bx|) = \frac 1{(2\pi\sigma)^{n/2}} \exp\left\{ - \frac{|\bx|^2}{2\sigma}\right\}.
\ee
Then
\[
K(|\bx|) = \displaystyle \frac{\int_{|\bx|^2}^{+\infty}\exp\left\{ - \frac{y}{2\sigma}\right\} \, dy}{2 \exp\left\{ - \frac{|\bx|^2}{2\sigma}\right\}} =\sigma.
\]
This suggests that Poincar\'e inequality for a Gaussian density is characterized by a constant weight. This is a very well-known result that easily follows from the logarithmic Sobolev inequality \cite{OV}. \\
A second important example is obtained when the isotropic density function is Cauchy--type density in $\R^n$ \cite{BL} 
\be\label{gg}
f_\infty(|\bx|) = \frac {C_\beta}{(1 +|\bx|^2)^\beta}, \quad \beta > \frac n2.
\ee
In \fer{gg} the constant $C_\beta$ is chosen to fix the mass equal to unity. In this case
\be\label{fat-k}
K(|\bx|) = \displaystyle \frac{\int_{|\bx|^2}^{+\infty}(1 +y)^{-\beta} \, dy}{2 (1 +|\bx|^2)^{-\beta}} =\frac 1{\beta-1}(1 +|\bx|^2) .
\ee
Hence, one can argue that for the Cauchy--type density  \fer{gg} the right weight in the Poincar\'e inequality \fer{Poi}  is given by $w(|\bx|) = \frac c{\beta-1}(1 +|\bx|^2) $. This result has in fact been proven by Bobkov and Ledoux in \cite{BL}, unlikely without the explicit evaluation of the constant $c$. \\
A third example is concerned with the exponential--type density in $\R^n$
\be\label{exp}
f_\infty(|\bx|) = \frac{\beta^n \Gamma\left(\frac n2\right)}{2\pi^{n/2 }\Gamma(n)} \exp\{ -\beta|\bx|\}, \quad \beta >0.
\ee
The value of $K$ is easily found to be
\be\label{exp-k}
K(|\bx|) = \displaystyle \frac{\int_{|\bx|^2}^{+\infty}e^{-\beta\sqrt y} \, dy}{2 e^{-\beta|\bx|}} =\frac 1{\beta^2}(1 +\beta |\bx|) .
\ee
These examples allow us to conjecture that, for a given isotropic probability density defined in $\R^n$, the value of the weight function $K(|\bx|)$ is an increasing function which grows at infinity like $|\bx|^p$, with $0\le p \le 2$, where the value 0 is attained in correspondence to the Gaussian density, while the value $2$ is attained in correspondence to the generalized Gaussians, which exhibit a polynomial decreasing at infinity.

A further interesting case is furnished by the class of Barenblatt solutions of the porous medium equation \cite{Vaz}, which, for  given constants $a>0$ and $p>1$ are defined as
\be\label{Bar-n}
f_\infty(|\bx|) = C_{a,p} (a^2 - |\bx|^2)^{\frac 1{p-1}},  
\ee
where  the constant $C_{a,p}$ is chosen to fix the mass equal to unity. In this case, $\MI$ is the ball of radius $a$, and,  for $ |\bx| \le a$ 
\be\label{Bar-k}
K(|\bx|) = \displaystyle \frac{\int_{|\bx|^2}^{a^2}(a^2 -y)^{1/(p-1)} \, dy}{2 (a^2 - |\bx|^2)^{1/(p-1)}} =\frac{p-1}{2p}(a^2 -|\bx|^2).
\ee
In what follows, we will show that a Fokker--Planck description of steady states allows us to obtain a proof of the validity of Poincar\'e inequalities with weight in many situations. The method of proof is an extension of the one-dimensional proof introduced in \cite{FPTT} and in the case of the Gaussian and  Cauchy--type densities the weight we found is exactly the one conjectured   here above.
Unfortunately in other examples we were not able to reach the expected sharp weight. Nevertheless, we were able to find the correct weight but only at the prize of restricting the domain of test functions in the Poincaré--type inequality (Proposition \ref{C}).


\section{A first Poincar\'e--type result for isotropic densities }\label{sec:poincare}

To start with, we prove that a Poincar\'e--type inequality  can be recovered by suitably generalizing the one-dimensional Poincar\'e inequalities obtained in \cite{FPTT21}. 
\subsection{A Poincar\'e--type inequality for product probability densities}\label{sec:prod}
The results in \cite{FPTT21} refer to  a class  of Fokker--Planck type equations with variable coefficients of diffusion and drift in the form
\begin{equation}\label{FP-gen}
 \frac{\partial f(x,t)}{\partial t} =   \frac{\partial^2 }{\partial x^2}
 \left(P(x)f(x,t)\right )+ 
\, \frac{\partial}{\partial x}\left( Q(x) f(x,t)\right),
 \end{equation}
where $x \in \MI =(a,b) \subseteq \R$, with $a<b$, coupled with no-flux boundary conditions
\be\label{bc}
\left. \frac{\partial }{\partial x}
 \left(P(x)f(x,t)\right )+ \, Q(x) f(x,t)\right|_{x= \{a, b\}} = 0,
\ee
which ensure mass conservation. In equation \fer{FP-gen} the coefficients of the diffusion $P(x)$  and the drift $Q(x)$ are smooth  functions, and $P(x) > 0$. In addition, $P(x)$ and $Q(x)$ are such that, for any $\alpha\in \MI$, the steady state 
 \be\label{eq-11}
 f_\infty(x) = \frac{C_\alpha}{P(x)} \exp \left\{  - \int_\alpha^x\, \frac{Q(y)}{P(y)}\, dy\right\}
 \ee
 is a probability density supported in $\MI$ for a given value of the constant $C_\alpha$.
Clearly, the steady state of the Fokker--Planck \fer{FP-gen}  solves the first order differential equation
\begin{equation}\label{ste1}
   \frac{\partial }{\partial x}
 \left(P(x)f_\infty(x)\right )+ 
\,  Q(x) f_\infty(x) = 0.
 \end{equation}
The case 
\be\label{facile}
Q(x) = x-m, \quad m \in (a,b)
\ee
has  already been considered and studied in \cite{FPTT}.  In \cite{FPTT} it was proven that, 
if $X$ is a random variable distributed with density $f_\infty(x)$, $x \in \MI \subseteq \R$, and $f_\infty$ satisfies the differential equation
\be\label{staz-22}
\frac{\partial }{\partial x}\left(P(x)  f_\infty(x) \right) + (x -m)\,f_\infty(x) = 0, \quad x\in \MI, 
 \ee
then for any smooth function $\phi$  defined  on $\MI$ such that $\phi(X)$ has finite variance, it holds
\be\label{ch-lin}
 Var[\phi(X)] \le E\left\{P(X)[\phi'(X)]^2\right\}.
\ee
Let us now reverse our point of view and consider
$f: (a,b) \to \R, (a,b) \subseteq \R $ be a probability density with finite mean value, say $m$, i.e. satisfying the condition
$\int_a^b x f(x) dx =m$, with $a<m<b$. It is immediate to verify that $f$ satisfies \fer{staz-22}, provided that 
\begin{equation}\label{Ca}
P(x)=\left\{
	\begin{array}{ll}
\dfrac{ \int_{a}^{x} (m-y)f(y)dy}{f(x)}, & a<x\leq m 
	 \\
	\dfrac{ \int_{x}^{b} (y-m)f(y)dy}{f(x)}, & m<x< b.
	\\ 
  	\end{array}
\right.
\end{equation}
We remark that $P$ is continuous on $(a,b)$. 
Therefore, if $X$ is a random variable distributed with density $f$, $P$, as given by \fer{Ca}, is the correct weight to satisfy inequality \fer{ch-lin}. \\
The results in \cite{FPTT} were recently extended in \cite{FPTT21} to cover the case of a more general drift function $Q$ as follows.
\begin{thm}\cite{FPTT21}\label{MJM}
Let $X$ be a random variable distributed with density $f_\infty(x)$, $x \in \MI =(i_-, i_+)\subseteq \R$. Let us suppose moreover that $f_\infty$ satisfies \fer{ste1}, where $Q(x)$ is a $C^1$ function, satisfying $Q'(x) >0$ on $\MI$, such that \be\label{bordi}
 \lim_{x\to i_-} Q(x)<0,\quad \lim_{x\to i_+} Q(x) >0.
 \ee
 Let $w(x)$ be defined by  
 \be\label{mag}
w(x) = \frac {P(x)}{Q'(x)},\quad x\in \MI.
 \ee
 Then,  for any  smooth function $\phi$  on $\MI$, and $\phi(X)$ with finite variance, it holds  
\[
 Var[\phi(X)] \le E\left\{w (X)[\phi'(X)]^2\right\},
\]
that is 
 \be\label{ch-gen}
 \int_\MI \phi^2(x) f_\infty(x) \, d x -\left(  \int_\MI \phi(x) f_\infty(x) \, d x\right )^2 \leq   \int_\MI w(x) \left(\phi'(x)\right )^2 f_\infty(x) \, d x.
 \ee
\end{thm}
\begin{rem} 
Of course for $f_\infty$ a probability density we have
\[
\int_\MI \phi^2(x) f_\infty(x) \, d x -\left(  \int_\MI \phi(x) f_\infty(x) \, d x\right )^2 =  \int_{\MI} \left(\phi(x) -  \int_{\MI} \phi(y) f_\infty(y) \, d y\right )^2 f_\infty(x) \, dx.
\]
\end{rem}
Note that inequality \fer{ch-gen} contains \fer{ch-lin} as a particular case, since $Q(x) = x-m$ implies $Q'(x) = 1$, and, in view of \fer{mag} $w(x) = P(x)$. \\
Having in mind the one-dimensional result \fer{ch-gen}, in the special case \eqref{ch-lin},
let us now consider a probability density  
 $f(\bx)$, with $\bx = (x_1,\dots, x_n) \in \R^n$ that is defined by the product of one-dimensional probability density functions, i.e.
\be\label{prod1} 
f(\bx) = \prod_{i=1}^n f_i(x_i).
\ee
In \fer{prod1}  each $f_i(x_i)$ is a probability density supported on $(a_i, b_i) \subseteq \R,$ with  $-\infty \leq a_i < b_i \leq +\infty, i=1,\dots, n$ such that the momentum
$$m_i = \int_{a_i}^{b_i} x_i f_i(x_i) \,d x_i
$$
is finite.
Any density $f_i$ is a solution to the following differential equation:
\begin{equation*}
\partial_{x_i}(P_i(x_i) f_i(x_i)) + (x_i -m_i)f_i(x_i)=0,
\end{equation*}
where $P_i(x_i)$ is defined as in \eqref{Ca}.
If $X_i$ is a random variable distributed with density $f_i(x_i)$, then $X_i$  satisfies the Poincar\'e inequality \fer{ch-gen} with weight function $P_i(x_i)$. \\
Let us denote by $A_n$ the set in $\R^n$ defined by
$$
A_n = \prod_{i=1}^n (a_i,b_i).
$$
For a product probability density we can prove the following  result.
\begin{thm}\label{Ch_Ind}
Let us suppose that a smooth scalar function $\phi(\bx)$, $\bx \in A_n$, has finite va\-riance  with respect to $f(\bx) = \prod_{i=1}^n f_i(x_i).$
If a random variable $X_i$ distributed with density $f_i$ satisfies a weighted one-dimensional  Poincar\'e inequality \fer{Poi} on $(a_i,b_i)$, with weight $w_i$, then 
  it holds 
\begin{equation}\label{Ch-product}
 \int_{A_n} \left(\phi(\bx) -  \int_{A_n} \phi(\by) f(\by) \, d \by\right )^2 f(\bx) \, d\bx  \leq   
  \int_{A_n} \sum_{i=1}^n w_i(x_i)\left[\frac{ \partial \phi (\bx) }{\partial x_i}\right]^2 f(\bx)\, d \bx.
\end{equation}
\end{thm}
\begin{rem}
In view of the discussion above, every density $f_i$  is associated with a weighted Poincaré inequality with weight $P_i$ defined as in \eqref{Ca}.  Moreover, there could be other weighted Poincaré--type inequalities associated with  the same density, as we will see in the examples in Section \eqref{sec:examples}.
\end{rem}

\begin{proof}
We prove the statement by induction with respect to the dimension $n.$
For $n=1$ it holds by virtue of the Poincaré inequality with weight proven in \cite{FPTT}. Suppose now that \fer{Ch-product}  is true for $n-1.$ For every fixed $x_n \in (a_n,b_n)$, let $\tilde\bx =(x_1,\dots, x_{n-1}) \in A_{n-1}$, and define the function 
$\phi_{x_n}$  on $A_{n-1}$ as
\begin{equation}
\phi_{x_n}(\tilde\bx)= \phi (\tilde\bx, x_n)= \phi(\bx).
\end{equation}
Due to the induction hypothesis,  the theorem is true for the function $\phi_{x_n},$ since for almost all $x_n$   $Var(\phi_{x_n})< +\infty $. Therefore 
\begin{equation*}
\begin{aligned}
&\int_{A_{n-1}}
\left[ \phi_{x_n}(\tilde\bx)  - \int_{A_{n-1}} \phi_{x_n}( \tilde\by)\prod_{i=1}^{n-1} f_i(y_i)d y_1\dots dy_{n-1}\right]^2 \prod_{i=1}^{n-1} f_i(x_i)dx_1\dots dx_{n-1} \leq \\
&\sum_{i=1}^{n-1}\int_{A_{n-1}}w_i(x_i)\left[\frac{ \partial \phi_{x_n} }{\partial x_i}\right]^2 \prod_{i=1}^{n-1}f_i(x_i)dx_1\dots dx_{n-1}, 
\end{aligned}
\end{equation*} 
where $\tilde\by =(y_1,\dots, y_{n-1}) \in A_{n-1}$.
By definition
\begin{equation*}
\begin{aligned}
Var(\phi(\bX)) = &\int_{A_n} \left(\phi(\bx) -  \int_{A_n} \phi(\by) f(\by) \, d \by\right )^2 f(\bx) d\bx  =\\
 &\int_{a_n}^{b_n}\int_{A_{n-1}} \phi^2(\bx) \prod_{i=1}^{n-1}f_i(x_i)\, f_n(x_n)dx_1 \dots dx_n -\\
 &\left (\int_{a_n}^{b_n} \int_{A_{n-1}} \phi(\bx)\prod_{i=1}^{n-1} f_i(x_i)\, f_n(x_n) d x_1 \dots d x_n \right)^2.
 \end{aligned}
\end{equation*}
Let us rewrite the previous equality 
by adding and subtracting the  quantity 
\[
 \int_{a_n}^{b_n}\left(\int_{A_{n-1}} \phi_{x_n}( \tilde\by)\prod_{i=1}^{n-1} f_i(y_i)d y_1\dots dy_{n-1}\right)^2f_n(x_n)\, dx_n.
\]
We obtain
\begin{equation*} 
\begin{aligned}
&\int_{a_n}^{b_n}\int_{A_{n-1}} \phi^2(\bx) \prod_{i=1}^{n-1}f_i(x_i)  f_n(x_n)dx_1\dots dx_n \\
&-\int_{a_n}^{b_n} \left( \int_{A_{n-1}}\phi_{x_n}( \tilde\by)
\prod_{i=1}^{n-1}f_i(y_i) dy_1 \dots dy_{n-1} \right)^2 f_n(x_n) dx_n \\
&+\int_{a_n}^{b_n}\left( \int_{A_{n-1}} \phi_{x_n}( \tilde\by)
\prod_{i=1}^{n-1} f_i(y_i) d y_1 \dots dy_{n-1} \right)^2 f_n(x_n) dx_n-\\
&\left[\int_{a_n}^{b_n} \left(\int_{A_{n-1}} \phi(\tilde\by,x_n)\prod_{i=1}^{n-1} f_i(y_i)d y_1 \dots dy_{n-1}\right ) f_n(x_n)dx_n \right]^2 .\\
 \end{aligned}
\end{equation*}
Thanks to  the induction hypothesis, from the first two terms we get the estimate
\begin{equation} \label{primi}
\begin{aligned}
&\int_{a_n}^{b_n}f_n(x_n) dx_n \cdot\\
&\cdot \left[ \int_{A_{n-1}}\phi^2_{x_n}(\tilde\bx) \prod_{i=1}^{n-1}f_i(x_i) dx_1 \dots dx_{n-1}- 
\left(\int_{A_{n-1}}\phi_{x_n}(\tilde\by) 
\prod_{i=1}^{n-1}f_i(y_i)dy_1\dots dy_{n-1}\right)^2 \right]  \leq\\
&\int_{a_n}^{b_n}f_n(x_n) dx_n \int_{A_{n-1}}\sum_{i=1}^{n-1} w_i(x_i)\left( \frac{\partial \phi}{\partial x_i}  \right)^2\ \prod_{i=1}^{n-1}f_i(x_i)dx_1 \dots dx_{n-1} =\\
&\int_{A_n} \left[\sum_{i=1}^{n-1}w_i(x_i)\left( \frac{\partial \phi}{\partial x_i}  \right)^2\right]\prod_{i=1}^{n}f_i(x_i) dx_1 \dots dx_n.\\
\end{aligned}
\end{equation}
Let us take now into account the difference of the last two terms, i.e.
\begin{equation*}
\begin{aligned}
&\int_{a_n}^{b_n} \left (\int_{A_{n-1}} \phi(\tilde\by,x_n)\prod_{i=1}^{n-1} f_i(y_i) d y_1 \dots dy_{n-1} \right)^2 f_n(x_n) dx_n -\\
&\left[ \int_{a_n}^{b_n}  \left (\int_{A_{n-1}} \phi(\tilde\by,x_n)\prod_{i=1}^{n-1} f_i(y_i) d y_1 \dots dy_{n-1}\right )f_n(x_n)dx_n\right]^2. 
 \end{aligned}
\end{equation*}
If we define
\[
{\tilde \phi} (x_n) =\int_{A_{n-1}} \phi(\tilde\by,x_n) \prod_{i=1}^{n-1} f_i(y_i) d y_1 \dots d y_{n-1},
\]
 we obtain that the previous difference is nothing but the  variance of the one-dimensional function $\tilde \phi (x_n)$. Indeed
\begin{equation*}
\begin{aligned}
Var(\tilde\phi(\bX_n)) =& \int_{a_n}^{b_n}\left( {\tilde \phi} (x_n) -\int_{a_n}^{b_n}{\tilde \phi }(y_n) f_n(y_n) dy_n \right)^2 f_n(x_n)dx_n = \\& \int_{a_n}^{b_n}
{\tilde \phi}^2 (x_n) f_n(x_n) dx_n  -\left (\int_{a_n}^{b_n}{\tilde \phi }(x_n) f_n(x_n) dx_n\right)^2.
\end{aligned}
\end{equation*}
Since the one-dimensional function $\tilde \phi(x_n)$  satisfies the hypotheses 
of the theorem about Poincaré inequality proven in \cite{FPTT} 
we get
\begin{equation*}
\begin{aligned}
& \int_{a_n}^{b_n}\left( {\tilde \phi} (x_n) -\int_{a_n}^{b_n}{\tilde \phi }(y_n) f_n(y_n) dy_n \right)^2 f_n(x_n)dx_n \leq \\
&\int_{a_n}^{b_n} w_n(x_n)\left(\frac{\partial \tilde \phi}{\partial x_n}\right)^2 f_n(x_n) d x_n.
 \end{aligned}
\end{equation*}
By Jensen inequality we have
\[
\begin{aligned}
\left(\frac{\partial \tilde \phi}{\partial x_n}\right)^2  &=	\left( \int_{A_{n-1}}  \frac {\partial }{\partial x_n} \phi(\tilde\by,x_n) \prod_{i=1}^{n-1} f_i(y_i) d y_1 \dots d y_{n-1}\right )^2	\\
&\leq \int_{A_{n-1}} \left( \frac {\partial  \phi }{\partial x_n}\right )^2 \ \prod_{i=1}^{n-1} f_i(y_i) d y_1 \dots d y_{n-1}.
\end{aligned}
\]
Finally we obtain
\be\label{secondi}
\begin{split}
 \int_{a_n}^{b_n}\left( {\tilde \phi} (x_n) -\int_{a_n}^{b_n}{\tilde \phi }(y_n) f_n(y_n) dy_n \right)^2 f_n(x_n)dx_n\\ \leq \int_{A_{n}}  w_n(x_n)\left(\frac{\partial\phi}{\partial x_n}\right)^2 \prod_{i=1}^{n} f_i(x_i) dx_1\dots dx_n.
 \end{split}
 \ee
 Gathering the two inequalities \eqref{primi} and \eqref{secondi} we complete the proof.
\end{proof}
\subsection{Poincar\'e--type inequalities for radially symmetric densities in $\R^n$.}\label{sec:chernoff}
The result of the previous section allows us to  prove a first Poincar\'e--type inequality with an explicit weight for probability densities in $\R^n$ which are radially symmetric. The reason is that in this case we can express the radially symmetric density as the product of $n$ one-dimensional densities in a system of hyperspherical coordinates. \\
 Before proving this result, let us briefly  recall  the relationship between the gradient in the Cartesian and in the hyperspherical coordinates. Let us denote by $\bx=(x_1,\dots, x_n)$  the Cartesian coordinates and  by $(\rho, \Theta) = (\rho, \theta_1, \dots, \theta_{n-1})$  the hyperspherical coordinates where $\rho =|\bx| \in(0,+\infty)$, $\theta_i\in (0,\pi)$, for any $i=1,\dots, n-2$ and  $\theta_{n-1}\in (0,2\pi)$. Then the Cartesian coordinates can be expressed in terms of the hyperspherical ones. We have
\begin{equation*}
\left\{
\begin{aligned}
&x_n=\rho \sin \theta_1 \dots \sin \theta_{n-3}\sin \theta_{n-2}\sin \theta_{n-1}\\
&x_{n-1}=\rho \sin \theta_1 \dots \sin \theta_{n-3}\sin \theta_{n-2}\cos \theta_{n-1}\\
&\vdots\\
&x_{3}=\rho \sin \theta_1 \sin\theta_2\cos \theta_{3}\\
&x_2=\rho \sin \theta_1  \cos \theta_{2}\\
&x_1=\rho \cos \theta_{1}.\\
\end{aligned}
\right.
\end{equation*}
Let $f:\R^n\to \R$, $f\in C^1$, $n\geq 2$ and let us denote 
\begin{equation}\label{R}
{\hat f}(\rho,\Theta)=f(x_1(\rho,\Theta), \dots, x_n(\rho,\Theta)).
\end{equation}
Then, the gradient with respect to the hyperspherical  coordinates is given by
\be\label{Gr}
	(\nabla f)({\bf x}(\rho,\Theta))=\frac{\partial{\hat f}}{\partial \rho}{e_\rho}+\frac 1\rho\frac{\partial{\hat f}}{\partial \theta_1}{ e_{\theta_1}}+\frac 1{\rho\sin \theta_1}\frac{\partial{\hat f}}{\partial \theta_2}{ e_{\theta_2}}+\dots+
	\frac 1{\rho\sin{\theta_1}\dots\sin{\theta_{n-2}}}\frac{\partial{\hat f}}{\partial \theta_{n-1}}{ e_{\theta_{n-1}}}
	\ee
	where $
	(e_\rho, e_{\theta_1},
	e_{\theta_2}, \dots,e_{\theta_{n-1}})
	$ is an orthonormal basis  of vectors.\\
The main result of this section is contained in the following theorem.
\begin{thm}\label{Cher}
	Let $\bX$ be a random variable distributed with a radially symmetric density 
	$f_\infty(|\bx|)$, $\bx \in \R^n$ and let $f(\rho)$, $\rho \in \R_+$, denote the one dimensional probability density function 
	\begin{equation}\label{H1}
	f(\rho) =\sigma_n \rho^{n-1} f_\infty(\rho),
	\end{equation} 
	where   $\sigma_n={2 \pi^{\frac n2}}/{\Gamma (\frac n2)}$ is the measure of the hyperspherical surface of unit radius  in $\R^n$. Then, if the non-negative random variable $X$ distributed with probability density $f(\rho)$ satisfies the one-dimensional  Poincar\'e inequality \fer{Poi} on ${\it I}=(0,+\infty)$, with weight $w(\rho)$,  for any smooth scalar function $\phi(\bx)$   such that $Var\left(\phi(\bX) \right)$ is finite it holds 
	\begin{equation}\label{Ch-gen}
	\int_{\R^n} \phi^2(\bx)f_\infty(|\bx|) \, d\bx -  \left(\int_{\R^n} \phi(\by) f_\infty(|\by|) \, d \by\right )^2 \leq  \int_{\R^n} W_*(|x|) \left|\nabla\phi(\bx)\right |^2 f_\infty(|\bx|) \, d \bx
	\end{equation}
	where 
	\[ 
	W_*(|x|)= \max \left\{w(|\bx|),\frac{\pi^2}2 |\bx|^2 \right\}.
    	\]

\end{thm}

\begin{proof} For the sake of simplicity we will prove the result in $\R^3$. Then the general proof will follow along the same lines.
The proof  is an easy application of Theorem \ref{Ch_Ind}. 	Indeed, let us evaluate  the left-hand side of  inequality \eqref{Ch-gen} by resorting to spherical coordinates. If $\bX$ is a random vector distributed with density $f_\infty(|\bx|)$, with $\bx \in \R^3$, then, if $\phi$ has finite variance, 
we write
	\begin{equation*}
	\begin{aligned}
	Var[\phi(\bX)]= \int_{\R^3} \phi^2(\bx)f_\infty(\bx) \, d\bx& -  \left(\int_{\R^3} \phi(\by) f_\infty(\by) \, d \by\right )^2=\\ =
	\int_0^{\infty}\int_0^\pi\int_0^{2\pi}\phi^2(\rho,\theta_1,\theta_2)& f_\infty(\rho)\rho^2\sin\theta_1 \, d\theta_1\, d\theta_2\, d\rho \\&-\left(\int_0^{\infty}\int_0^\pi\int_0^{2\pi}\phi(\rho,\theta_1,\theta_2)  f_\infty(\rho)\rho^2\sin\theta_1 \, d\theta_1\, d\theta_2\, d\rho \right)^2\\=
	\int_0^{\infty}\int_0^\pi\int_0^{2\pi}\phi^2(\rho,\theta_1,\theta_2)& f(\rho)f_1(\theta_1)f_2(\theta_2) \, d\theta_1\, d\theta_2\, d\rho \\
	&-\Big(\int_0^{\infty}\int_0^\pi\int_0^{2\pi}\phi(\rho,\theta_1,\theta_2) f(\rho)f_1(\theta_1)f_2(\theta_2) \, d\theta_1\, d\theta_2\, d\rho \Big)^2,
	\end{aligned}
	\end{equation*}
	where $f$, $f_1$ and $f_2$ are \emph{probability densities} defined by
	\begin{equation}\label{f}
	f(\rho)=4\pi\rho^2 f_\infty(\rho),\quad \rho \in (0,\infty), 
	\end{equation}
	and
	\begin{equation}\label{f1}
	 f_1(\theta_1)=\frac{\sin \theta_1}{2},  \quad \theta_1\in(0,\pi); \qquad f_2(\theta_2)=\frac1{2\pi}, \quad \theta_2\in(0,2\pi).
	\end{equation}
	By assumption, the probability density $f(\rho)$ satisfies \fer{ch-gen} with a weight function 
	$w(\rho) $. \\
	Furthermore the densities 
	of $f_1$ and $f_2$ have finite mean values
	$$
	m_1=\int_0^\pi\theta_1\frac{\sin\theta_1}2d\theta_1=\frac \pi2 \ \ \    \hbox{and }\ m_2=\int_0^{2\pi}\frac{\theta_2}{2\pi}d\theta_2=\pi
	$$
	and they respectively solve 
	on $(a_i,b_i)$ with $(a_1,b_1)=(0,\pi)$ and  $(a_2,b_2)=(0,2\pi)$ the differential equations
	\begin{equation}\label {K1}
	\frac{\partial}{\partial\theta_i}\left[P_i(\theta_i)f_i(\theta_i)\right]+(\theta_i-m_i)f_i(\theta_i)=0
	\end{equation}
	where $P_i(\theta_i)$ are defined as in \eqref{Ca}. \\
	Since the function $\sin \lambda$ is increasing in the interval $(0, \pi/2)$, and decreasing  in the interval $( \pi/2, \pi)$ we have
	\begin{equation}\label{teta1}
\dfrac{ \int_{0}^{\theta_1} \left(\frac \pi 2- \lambda\right)\sin\lambda\, d\lambda}{\sin\theta_1} \le { \int_{0}^{\theta_1} \left(\frac \pi 2- \lambda\right)\, d\lambda} \le \frac{\pi^2}{8},  \quad  0< \theta_1\leq\frac \pi 2.
	\end{equation}
	Likewise	
\begin{equation}\label{teta11}
	\dfrac{ \int_{\theta_1}^\pi \left( \lambda - \frac \pi 2\right)\sin\lambda\, d\lambda}{\sin\theta_1} \le { \int_{\theta_1}^\pi \left( \lambda- \frac \pi 2\right)\, d\lambda} \le \frac{\pi^2}{8},  \quad  \frac \pi 2< \theta_1\leq \pi .
	\end{equation}
	Consequently 
	\be\label{k1}
	P_1(\theta_1) \le \frac{\pi^2}{8}, \quad 0 \le \theta_1 \le \pi.
	\ee
	Also $$
	P_2(\theta_2)=\pi\theta_2-\frac{\theta_2^2}{2} \le \frac{\pi^2}{2}, \quad 0 \le \theta_2 \le 2\pi.
	$$
	Applying  Theorem \ref{Ch_Ind} to the density $f(\rho) f_1(\theta_1) f_2(\theta_2)$ on $(0,+\infty)\times (0,\pi)\times(0,2\pi)$, we conclude that 	
	\begin{equation}\label{Chern-sferiche}
	Var[\phi(\bX)] \leq\int_0^\infty\int_0^\pi\int_0^{2\pi}\Bigg[w(\rho) \left|\frac{\partial \phi}{\partial\rho}\right|^2+ \frac{\pi^2}{8} \left|\frac{\partial \phi}{\partial\theta_1}\right|^2
	+ \frac{\pi^2}{2}\left|\frac{\partial \phi}{\partial\theta_2}\right|^2\Bigg]f(\rho)f_1(\theta_1)f_2(\theta_2) \, d\rho\, d\theta_1\, d\theta_2,
	\end{equation}
	with $ \left|\frac{\partial \phi}{\partial\rho}\right|^2 \leq |\nabla \phi|^2 $ and 
	$ \left|\frac{\partial \phi}{\partial\theta_i}\right|^2  \leq   \rho^2 |\nabla \phi|^2, \,  i=1,2.$
		Reverting back to Cartesian coordinates we finally obtain 
	\be
	I \leq \int_{\R^3}\max \left\{w(|\bx|),  \frac{\pi^2}{2} |\bx|^2\right\} |\nabla \phi|^2f_\infty(|\bx|)\, d\bx.
	\ee	
	\end{proof}
	\begin{rem} It is important to remark that in the general case the upper bound on the coefficients $P_i(\theta_i)$, $i =1,\dots, n-1$ does not change, since the probability density function leading to $P_i(\theta_i)$, $0 \le \theta_i \le \pi$ is given by
	\be
	f_i(\theta_i) = \frac{\sin^i \theta_i}{\int_0^\pi \sin^i \lambda\, d\lambda }, \quad i = 1, 2,\dots, n-2
	\ee
	and we can always  reduce to  inequalities like \fer{teta1} and \fer{teta11}. 
		\end{rem}
	\begin{rem} The previous result also holds  for radially symmetric density 
			$f_\infty(|\bx|)$ defined on a sphere  $B_R=\{\bx \in \R^n:\ |x|\leq R\} $, for some $R>0$. 
			\end{rem}
\begin{rem} As we shall see in Section \ref{sec:examples}, there exists in general a whole family of Fokker--Planck equations which share the same stationary state $f(\rho) = \sigma_n \rho^{n-1} f_\infty(\rho)$. The possibility of choosing which equation is more suitable in order to optimize the weight in the Poincaré inequality is one of the main feature of this line of research.
\end{rem}

%
%

\subsection{Examples.}\label{sec:examples}

In this section we collect some examples  in which the previous results  can be fruitfully applied to get Poincar\'e--type inequality with weight \fer{Poi}. In particular, we will study in details the probability density functions considered in Section \ref{sec:Peso}. \\
The analysis of Section \ref{sec:chernoff} leads to the conclusion that the major difficulty in getting a Poincar\'e inequality  for a radially symmetric density is the determination of the weight function $w(\rho)$ of the radial density \eqref{H1}.  
\subsubsection{The case of a Gaussian density function}\label{sec:gaussian} 

Given  a $n\times n$ real symmetric, positive definite matrix $V$ and a vector ${\bf {u}} =(u_1,\dots,u_n),$ let us consider
the multivariate Gaussian density $G$ 
$$
G({\bf x})=\frac{1}{(2\pi)^{n/2}}|V|^{-1/2}{\rm exp}\left\{-\frac 12 ({\bf x}-{\bf {u}})^t\ V^{-1}\ ({\bf x}-{\bf {u}})\right\}
$$
which has $V$ as covariance matrix with  $|V|= det (V) \neq 0 $ and ${\bf {u}}$ as mean vector. 

\begin{thm}
	For any bounded, smooth function $\phi: \R^n \to \R,$ such that the variance of $\phi$  with respect to $G(\bf x)$ is finite,  it holds 
	\begin{equation}\label{in}
	\int_{\R^n} \left(\phi ({\bf x}) -  \int_{\R^n} \phi({\bf y}) G({\bf y}) \, d {\bf  y}\right )^2 G({\bf x}) d{\bf x}  \leq   
	\int_{\R^n} 
	\left(\max_{i=1, \dots n}\lambda_i\right)\left|\nabla\phi ({\bf x}) \right|^2 G({\bf x}) d {\bf x},
	\end{equation}
	where $\lambda_1, \dots, \lambda_n$ are the eigenvalues of the symmetric and  positive definite real matrix $V$.
\end{thm}

This is a classical result and a proof can be found for example in \cite{Chen}. We detail it below for the reader's convenience.

\begin{proof}
	It is well known that 
	given  a real symmetric, positive definite,  matrix $V, $ there exists a non singular matrix $H$ such that if 
	$$
	{\bf x}={\bf {u}}+H{\bf{x}^*},
	$$
	then the positive definite form in the exponent of the Gaussian becomes
	$$
	({\bf x}-{\bf {u}})^t\ V^{-1}\ ({\bf x} -{\bf {u}}) = \sum_{i=1}^n ({x_i}^*)^2.
	$$
	The matrix $H$ is given by
	$
	H=(O^t ) {\sqrt D}.
	$
	where $O$ is the orthogonal matrix which diagonalizes the real symmetric matrix $V$ and 
	$$
	D=OVO^t=\begin{pmatrix}
	\lambda_1&0&\cdots&0\\
	0&\lambda_2&\cdots&0\\
	&&\cdots&\\
	0&0&\cdots&\lambda_n\\
	\end{pmatrix}
	$$
	with  $\lambda_i$  eigenvalues of V.

	Let us consider the  previous change of variables
	$$
	{\mathbf x}={\bf{u}}+H{\bf {x}^*}.
	$$
	As $H(H^t)=V, $ we get $|H|=|V|^{\frac{1}{2}}.$ 
	We get
	\begin{equation*}
	\begin{aligned}
	&\int_{\R^n} \left( \phi({\bf x}) -  \int_{\R^n} \phi({\bf y}) G({\bf y}) \, d {\bf y}\right)^2 G({\bf x}) d{\bf x }=\\
	&\int_{\R^n} \left(\phi( {\mathbf u}+ H{\bf x}^*) 
	-  \int_{\R^n} \phi\left(
	{\mathbf u}+ H{\bf {y}}^*) G({\mathbf u}+ H{\bf y}^*\right)|H| \, d {\bf y}^ *
	\right)^2
	G({\mathbf u}+ H{\bf x}^*)|H| d{\bf x}^*= 
	\\
	&\int_{\R^n} \left(\Psi({\bf{x}}^*) -  \int_{\R^n} \Psi( {\bf {y}}^*) \frac 1{(\sqrt{2\pi})^n} \prod_{i=1}^{n} {\rm e}^{-\frac{(y_i^*)^2}{2}}\, d {\bf y}^*\right)^2\frac 1{(\sqrt{2\pi)}^n}
	\prod_{i=1}^{n}  {\rm e}^{-\frac{(x_i^*)^2}{2}} d {\bf x}^*,
	\end{aligned}
	\end{equation*}
	where, in the last line, we set $\Psi( {\bf {x}}^*) =\phi\left({\mathbf u}+ H{\bf {x}}^*\right ) .$
	As  any density 
	$$
	g(x_i)=\frac 1{\sqrt{2\pi}} {\rm e}^{-\frac{(x_i^*)^2}{2}}
	$$
	satisfies  the equation
	$$
	\partial_{x_i}( g(x_i))+x_ig(x_i)=0
	$$
	we can apply Theorem \ref{Ch_Ind} with $K_i=1, i=1,\dots ,n$ and obtain
	\begin{equation}\label{p3}
	\begin{aligned}
	&\int_{\R^n} \left(\Psi({\bf {x}}^*) -  \int_{\R^n} \Psi( {\bf {y}}^*)\frac 1{(\sqrt{2\pi)}^n}  \prod_{i=1}^{n} {\rm e}^{-\frac{(y_i^*)^2}{2}}\, d {\bf y}^*\right)^2\frac 1{(\sqrt{2\pi})^n}  \prod_{i=1}^{n} {\rm e}^{-\frac{(x_i^*)^2}{2}} d {\bf x}^*\\
	&\leq \int_{\R^n}\sum_{i=1}^n \left(\frac{\partial \Psi}{\partial{x_i^*}}(\bf x^*)\right)^2 \frac 1{(\sqrt{2\pi)}^n}  \prod_{i=1}^{n}  {\rm e}^{-\frac{(x_i^*)^2}{2}} d {\bf x}^*.
	\end{aligned}
	\end{equation}
	Now since $\Psi({\bf {x}^*})=\phi({\bf{u}}+H\bf x^*)$ we get
	$$
	\frac{\partial \Psi}{\partial{x_k^*}}=\sum_{i=1}^n\frac{\partial \phi}{\partial{x_i}} \frac{\partial x_i}{\partial{x_k^*}}=\sum_{i=1}^n(H^t)_{ki}\frac{\partial \phi}{\partial{x_i}} 
	$$
	and
	$$
	\nabla_{\bf {x^*}}\Psi = (H^t)\ \nabla_{\bf { x}}\phi
	$$ and
	\begin{equation}\label{p4}
	\left|\nabla_{\bf{ x^*}}\Psi\right|^2 =[(\nabla_{\bf{ x^*} }\Psi)^t] \nabla_{\bf{ x^*}}\Psi= 
	[\left(\nabla_{\bf{ x}}\phi\right)^t] H(H^t) \nabla_{\bf{ x}}\phi= 
	[(\nabla_{\bf{ x}}\phi)^t]V\nabla_{\bf{ x^*}}\phi,
	\end{equation}
	as $H(H^t)=V.$
	On the other hand, as $V= (O^t )D O,$ we get
	$$
	[(\nabla_{\bf{x}}\phi)^t] V [(\nabla_{\bf{ x}}\phi)] =({\bf {r}}^t )D {\bf{ r}},$$ where ${\bf{ r}} = 
	(r_1,\dots,r_n):= O\nabla_{\bf{x}}\phi.
	$ 
	As $O$ is an orthogonal matrix,  we get 
	\begin{equation}
	\left|\nabla_{\bf{ x^*}}\Psi\right|^2 = \sum_{i=1}^n \lambda_i r_i^2 \leq (\max_{i=1,\dots,n}\lambda_i)\sum_{i=1}^n r_i^2= (\max_{i=1,\dots, n}\lambda_i)
	|\nabla_{\bf{ x}}\phi |^2.
	\end{equation}
	Finally, thanks to \eqref{p4} and the change of variables ${\bf x}={\mathbf u}+H{\bf{ x*}}$,  for the last term in \eqref{p3} we get the estimate
	\begin{equation*}
	\begin{aligned}
	&\int_{\R^n}\sum_{i=1}^n \left(\frac{\partial \Psi}{\partial{x_i^*}}({\bf{ x^*}})\right)^2 \frac 1{(\sqrt{2\pi})^n} \prod_{i=1}^{n} {\rm e}^{-\frac{(x_i^*)^2}{2}} d {\bf{ x}^*}\\&= \int_{\R^n}\left|\nabla_{\bf{ x^*}}\Psi\right|^2\frac 1{(\sqrt{2\pi})^n}  \prod_{i=1}^{n}  {\rm e}^{-\frac{(x_i^*)^2}{2}} d {\bf{ x^*}}\\&\leq
	\int_{\R^n}\left(\max_{i=1, \dots, n} \lambda_i\right)\left|\nabla_{\bf {x}}\phi\right|^2 G({\bf x}) d {\bf{ x}}.
	\end{aligned}
	\end{equation*}

\end{proof}

In the examples which follow we can make use of the one--dimensional class of Fokker--Planck--type equations \fer{FP-gen} leading to the (optimal in the class) generalized Poincar\'e inequality \fer{ch-gen}. 

\subsubsection{Poincar\'e type inequalities for Cauchy-type densities}\label{sec:Cauchy}

The second example we consider is concerned with Cauchy densities,	
namely  probability densities with a polynomial decay at infinity given by
\begin{equation}\label{Cau}
f_\beta(|\bx|)=\frac{C_\beta}{(1+|\bx|^2)^\beta},\quad  \beta > \frac n2
\end{equation}
with $C_\beta$ satisfying
\[
\int_{\R^n}\frac{C_\beta}{(1+|\bx|^2)^\beta}d\bx=1.
\]
In this case, the radial density $f_\beta(\rho)$, $ \rho\in \R_+$  takes the form
\be\label{HH}
f_\beta (\rho) = \sigma_n \,\rho^{n-1}\frac{C_\beta}{(1+\rho^2)^\beta}, \quad \beta > \frac n2.
\ee
Following the strategy developed in \cite{FPTT21}, we observe that the Cauchy-type density $f_\beta$ defined in \fer{HH} can be characterized as the stationary state of a whole family of Fokker--Planck type equations, which depend on a positive parameter $\alpha$. For our purposes, we will assume that  $\alpha \in \left( 1/2, 1\right]$.  It can be easily checked that 
this family of Fokker--Planck type equations is given by
\[
\partial_t f_\beta = \partial^2_\rho \left((1+\rho^2)^\alpha f_\beta \right ) + \partial_\rho\left[ \left(\frac {2(\beta - \alpha) \rho}{(1+\rho^2)^{1-\alpha}} -\frac{n-1}\rho(1+\rho^2)^\alpha \right)f_\beta  \right], \quad \rho\in \R_+, t>0.
\]
Indeed, $f_\beta$ satisfies, for all $\rho\in\R_+$, the differential equation
\be\label{sta-b}
\partial_\rho \left((1+\rho^2)^\alpha f_\beta\right ) = - \left(\frac {2(\beta - \alpha) \rho}{(1+\rho^2)^{1-\alpha}} -\frac{n-1}\rho(1+\rho^2)^\alpha \right)f_\beta.
\ee
Equation \fer{sta-b} is of the type \fer{ste1}, with 
\be\label{pq}
P(\rho) = (1+\rho^2)^\alpha, \quad Q(\rho)= \frac {2(\beta - \alpha) \rho}{(1+\rho^2)^{1-\alpha}} -\frac{n-1}\rho(1+\rho^2)^\alpha, \quad \rho\in\R_+.
\ee
In the allowed range of the constant $\alpha$, the function $Q$ satisfies all the assumptions of Theorem \ref{MJM}.
Indeed,  for all $1/2 < \alpha \le 1$  and 
\be\label{cond}
\beta - \alpha > \frac{n-1}2,
\ee
the function $Q$ satisfy assumption \fer{bordi}, since
\[
\lim_{\rho\to 0^+} Q(\rho)=-\infty, \quad \lim_{\rho\to +\infty} Q(\rho)= +\infty.
\]
Condition \eqref{cond} holds true for $\alpha \le 1$ provided 
\be\label{cond2}
\beta > (n+1)/2.
\ee 
Moreover, $Q(\rho)$ is differentiable on $\R_+$ with
\[
Q'(\rho)= \frac {2(\beta -\alpha) (1+\rho^2(2\alpha -1))-2\alpha(n-1)(1+\rho^2)}{(1+\rho^2)^{2-\alpha}} + \frac{n-1}{\rho^2}(1+\rho^2)^\alpha, \quad  \rho\in \R_+.
\]
Let us observe that 
\[
\frac{n-1}{\rho^2}(1+\rho^2)^\alpha = \frac{(n-1)(1+\rho^2)^2 }{\rho^2(1+\rho^2)^{2-\alpha}} \ge \frac{(n-1)(1+\rho^2) }{(1+\rho^2)^{2-\alpha}}.
\]
Therefore
\[
Q'(\rho)\ge \frac {2(\beta -\alpha)(1+\rho^2(2\alpha -1))-(2\alpha-1) (n-1)(1+\rho^2)}{(1+\rho^2)^{2-\alpha}}.
\]
Moreover, since $\alpha \leq 1$ we have $2\alpha -1 \leq1$,   so that
\[
2(\beta -\alpha)(1+\rho^2(2\alpha -1)) \ge 2(2\alpha -1)(\beta -\alpha)(1+\rho^2).
\]
Hence
\[
Q'(\rho)\ge \frac {(2\alpha -1)\left[2(\beta -\alpha)- (n-1)\right]}{(1+\rho^2)^{1-\alpha}}.
\]
Hence, $Q'(\rho)$ is positive since condition \eqref{cond} holds true.
Under this condition $Q : \R \longrightarrow \R$ is a strictly monotone, smooth transformation.\\
Finally
\be\label{bound-c}
\frac{P(\rho)}{Q'(\rho)} \le   \frac 1{(2\alpha -1)\left[2(\beta -\alpha)- (n-1)\right]} (1+\rho^2),\quad \rho\in\R_+.
\ee
We can now look for the optimal value of the constant $(2\alpha -1)[2(\beta -\alpha)- (n-1)]$ under the constraint  $\alpha \in \left (1/2, 1\right ]$ and $\beta > (n+1)/2$. Let $\beta^* = \beta -(n-1)/2 > 1$.
The optimal value is obtained by maximizing the function
\[
h_\beta(\alpha)= (2\alpha -1)(\beta^*-\alpha).
\]
To this end, since
\[
h'_\beta(\alpha)= 2\beta^*-4\alpha+1
\]
we obtain
\[
h_\beta'(\alpha) \geq 0 \iff \alpha \leq \frac{\beta^*} 2+\frac 14.
\]
Hence
\[
\alpha_{\max}= \frac{\beta^*}2 +\frac 14, \quad{\rm if} \quad
\frac{\beta^*}2 +\frac 14 \leq 1
\]
while
\[
\alpha_{\max}=1, \quad{\rm if} \quad  \frac{\beta^*}2 +\frac 14>1.
\]
In conclusion we obtain
\be\label{w-c}
w(\beta, \rho) =\displaystyle
\begin{cases}
	& \frac 1{(\beta - n/2)^2} (1+\rho^2) \quad \frac n2 +\frac 12 < \beta < \frac n2 + 1, \\
	& \frac 1{2\beta -(n+1)}(1+\rho^2) \quad \beta \ge \frac n2 + 1.
\end{cases}
\ee
Therefore,  we apply Theorem \ref{Cher} and obtain
	\[
	W_*(|\bx|)  =   \max\left\{w(\beta, 0), \frac{\pi^2}2 \right \} (1 +|\bx|^2)
	\]
and $W_*(|\bx|) =c_{\beta,n} K(|x|)$ with $K(|x|)=\frac1{\beta-1}(1+|x|^2)$ as \eqref{fat-k}. 
This result is in perfect agreement with the analogous one proven by Bobkov and Ledoux in \cite{BL} without an explicit value of the constant in front of the weight. Moreover, in this case the weight is the one conjectured in Section \ref{sec:Peso}.

\vskip 3mm		

\subsubsection{Poincar\'e--type inequalities for exponentially decaying densities}\label{sec:exponential}

Let us now consider the case of the exponentially decaying  probability density function \fer{exp}, which is uniquely connected with the weight $K(|x|)=\frac{1}{\beta^2}(1+\beta |x|)$, as in \fer{exp-k}. In order  to explicitly determine the value of the weight $w(\rho)$, as given by \fer{mag}, it is essential to remark that in our case the radial function $f(\rho)$ in Theorem \ref{Cher} is  a Gamma density given by
\be\label{Gamma1}
f(\rho) = \frac{\beta^n}{\Gamma(n)} \rho^{n-1} e^{-\beta\rho}.
\ee
Poincar\'e-type inequalities with optimal weight for a large class of generalized Gamma densities have been recently obtained in \cite{To20}. As far as the density function \fer{Gamma1} is concerned,  one obtains $w(\rho) = \beta \rho$. Consequently, the weight in Theorem \ref{Cher} is given by
\[
W_*(|\bx|) = \max\{\beta |\bx|, \frac{\pi^2}2 |\bx|^2\}.
\]
We remark that, in this case, the application of Theorem \ref{Cher} leads to a non optimal weight, growing at infinity like $|\bx|^2$.

\subsubsection{Poincar\'e inequalities for Barenblatt type densities}\label{sec:Barenblatt}

Let us now consider the case of the class of radially symmetric Barenblatt type densities \fer{Bar-n}, which are characterized by a bounded spherical support, and by the weight function  $K(|x|)=\frac{p-1}{2p}(a^2-|x|^2)$ as in \fer{Bar-k}.\\
In what follows, for given $a >0$ and $p>1$, we consider the class of Barenblatt densities supported in the sphere $B_a=\{|\bx|\le a\} \subset \R^n$, $n >1$, in the form
\begin{equation}\label{Bar}
f_p(|\bx|)=C_{a,p}(a^2 -|\bx|^2)^{\frac1{p-1}}, \quad |\bx| \le a,
\end{equation}
where the constant  $C_{a,p}$ is chosen to fix the mass equal to unity.
For the sake of simplicity, we denote $\beta=\frac1{p-1}$ and with an abuse of notation we will write 
$f_\beta(|x|)=f_p(|x|)$ and $C_{a,\beta}= C_{a,p}$.
In this case, the radial density $g_\beta(\rho)$, $\rho \in \R_+$  takes the form
\be\label{HHH}
g_\beta (\rho) = \sigma_n \,\rho^{n-1}C_{a,{\beta}}{(a^2-\rho^2)^\beta}, \quad 0\le \rho \le a.
\ee
Following the strategy developed in Section \ref{sec:Cauchy}, we introduce a family of Fokker--Planck--type equations  that have, for a given $\beta >0$, the Barenblatt density $g_\beta$ defined in \fer{HHH}  as a stationary state. These equations depend on a positive parameter $\alpha$ that,  as in Section \ref{sec:Cauchy}, belongs to the interval  $ \left( 1/2, 1\right]$.  \\ 
It can be easily checked that 
this family of Fokker--Planck type equations is given by
\[
\partial_t g_\beta = \partial^2_\rho \left((a^2- \rho^2)^\alpha g_\beta \right ) + \partial_\rho\left[ \left(\frac {2( \alpha +\beta) \rho}{(a^2-\rho^2)^{1-\alpha}} -\frac{n-1}\rho(a^2- \rho^2)^\alpha \right)g_\beta  \right],
\]
for $  \rho\in \MI=(0,a)$ and $ t>0.$
Indeed, $g_\beta$ satisfies, for all $\rho\in\MI$, the differential equation
\be\label{sta-bar}
\partial_\rho \left((a^2 - \rho^2)^\alpha g_\beta\right ) = - \left(\frac {2( \alpha +\beta) \rho}{(a^2-\rho^2)^{1-\alpha}} -\frac{n-1}\rho(a^2- \rho^2)^\alpha \right)g_\beta.
\ee
Equation \fer{sta-bar} is of the type \fer{ste1}, with 
\be\label{pqb}
P(\rho) = (a^2 - \rho^2)^\alpha, \quad Q(\rho)= \frac {2( \alpha +\beta) \rho}{(a^2-\rho^2)^{1-\alpha}} -\frac{n-1}\rho(a^2- \rho^2)^\alpha , \quad \rho\in\MI.
\ee
In the allowed range of the constant $\alpha$, the function $Q$ satisfies all the assumptions of Theorem \ref{MJM}.
Indeed,  for all $1/2 < \alpha \le 1$  the function $Q$ satisfies assumption \fer{bordi}, since
\[
\lim_{\rho\to 0^+} Q(\rho)=-\infty, \quad \lim_{\rho\to a^-} Q(\rho)> 0.
\]
Moreover, $Q(\rho)$ is differentiable on $\MI$ with
\[
Q'(\rho)= \frac {2(\alpha +\beta) (a^2- \rho^2(2\alpha -1))}{(a^2-\rho^2)^{2-\alpha}} + \frac{n-1}{(a^2-\rho^2)^{1-\alpha}}\frac{a^2+(2\alpha -1)\rho^2}{\rho^2}, \quad  \rho\in\MI.
\]
Let us observe that, in the interval $ 1/2<  \alpha \le 1$ we have  $0 < 2\alpha -1\le 1$, so that 
\[
\frac {2(\alpha +\beta) (a^2- \rho^2(2\alpha -1))}{(a^2-\rho^2)^{2-\alpha}} \ge \frac {2(\alpha +\beta) (2\alpha -1)}{(a^2-\rho^2)^{1-\alpha}}.
\]
Moreover, for $\rho\in\MI$ we have
\[
\frac{n-1}{(a^2-\rho^2)^{1-\alpha}}\frac{a^2+(2\alpha-1)\rho^2}{\rho^2}\ge  \frac{2(n-1)(2\alpha-1)}{(a^2-\rho^2)^{1-\alpha}}.
\]
Finally,
\[
Q'(\rho)\ge \frac  {2(2\alpha-1)}{(a^2-\rho^2)^{1-\alpha}} (\alpha+\beta+n-1).
\]
This implies
\be\label{bound-bar}
\frac{P(\rho)}{Q'(\rho)} \le   \frac 1{2(2\alpha-1 )(\alpha+\beta+n-1) } (a^2-\rho^2),\quad \rho\in\MI.
\ee
We can now look for the optimal value of the constant  $2(2\alpha-1 )(\alpha+\beta+n-1) $ under the constraint  ${\alpha \in \left (1/2, 1\right ]}$ and $\beta > 0$. 
The optimal value is obtained by maximizing the function
\[
h_\beta(\alpha)= (2\alpha -1)(\alpha +\beta +n-1).
\]
To this end, since $\alpha> 1/2$
\[
h'_\beta(\alpha)= 2\beta+ 4\alpha-1 +2(n-1)>0
\]
and
\[
h_\beta(\alpha) \le 2(n+\beta).
\]
Going back to Definition \fer{mag} and recalling that $\beta(p)=\frac{1}{p-1},$  for $p>1$ we obtain
\be\label{w-bar}
w(p, \rho) \leq \frac 1{2(n+\beta(p))} (a^2- \rho^2) = \frac{p-1}{2(n(p-1)+1)}(a^2- \rho^2), 
\ee
and  in Theorem \ref{Cher}
\be\label{aaa}
W_*(|\bx|) = \max\left\{ \frac{p-1}{2(n(p-1)+1)} (a^2 -|\bx|^2) , \frac{\pi^2}2 |\bx|^2 \right\}.
\ee
Also in this case Theorem \ref{Cher} gives a non-optimal weight.
\vskip 3mm

\section{A refined Poincaré--type result for isotropic densities}\label{sec:refined}
The weight obtained in the Poincaré inequality \eqref{Ch-gen} grows at infinity at least as $|\bx|^2$ and this is due to the contribution of the angular components in the decomposition of the density.
This appears not to be sharp, if we consider the leading examples presented in  Section \ref{sec:Peso} and at the end of the previous section.  More specifically, in the case of the exponentially decaying density \eqref{exp} the decaying at infinity of the weight should be linear instead of quadratic and in the Barenblatt case \eqref{Bar} it should be vanishing at the boundary of the domain.
We are now presenting a refined Poincaré-type inequality with weight for a $n$-dimensional isotropic density. On one hand we sharpen the contribution of the weight for large scales but on the other we restrict the domain where it can be applied. 
\begin{prop}\label{C}
Let $\bX$ be a random vector distributed with density $f_\infty(|\bx|)$, with $\bx \in \R^n$  where the probability density function $f_\infty(|\bx|)$ satisfies the differential equation
\be\label{staz-2}
\nabla [K(|\bx|) f_\infty(|\bx|)] + \bx f_\infty(|\bx|) = 0.
 \ee
Moreover,  suppose  that there exists a real value $R>0$ such that, for any $r \ge R$, the weight function $K(r)$ satisfies the inequality 
 \be\label{b1}
(n-1) \frac{K(r)}{r^2} \le \frac 12. 
 \ee
Then, for any smooth function $\phi$, supported outside the ball $B_R=\{\bx \in \R^n: \ |\bx|\leq R \}$, and with finite variance with respect to $f_\infty(|\bx|)$, we have 
\[
\int_{\R^n} \phi^2(\bx)f_\infty(|\bx|) \, d\bx -  \left(\int_{\R^n} \phi(\by) f_\infty(|\by|) \, d \by\right )^2 \leq  2\int_{\R^n} K(|x|) \left|\nabla\phi(\bx)\right |^2 f_\infty(|\bx|) \, d \bx.
\]
\end{prop}
\begin{rem} We underline that in this proposition the weight we get agrees with the one conjectured in Section \ref{sec:Peso}. It would be interesting to obtain  the same result {\it without any restriction on the test function $\phi$}.
In the examples dealt with in Subsection \ref{sec:examples} the diffusion term $K$ always satisfies  condition \eqref{b1}.
\end{rem}

\begin{proof}
Let a smooth scalar function $\phi$ defined on $\R^n$ and such that $\phi({\bf x})=0$ for $|{\bf x}| \leq R$.
Since for any probability density function $f_\infty(|\bx|)$, we have the elementary inequality
 \[
\int_{{\R^n}} \phi ^2(\bx) f_\infty(|\bx|) \ d\bx -   \left( \int_{{\R^n}} \phi (\bx) f_\infty(|\bx|) \ d\bx \right )^2 \leq \int_{{\R^n}} \left(\phi (\bx)-\phi (0)\right )^2 f_\infty(|\bx|)\ d\bx,
\] 
we get
\[
\begin{split}
 Var[\phi (\bX)]  &=  \int_{\R^n}  \left(\phi (\bx) - \int_{\R^n} \phi  (\by) f_\infty(|\by|) \, d\by\right)^2 f_\infty(|\bx|)\, d\bx\\
 \leq & \int_{{\R^n}} (\phi (\bx)-\phi (0))^2 f_\infty(|\bx|)\ d\bx 
\end{split}
\]
and in hyperspherical coordinates
\[
 \int_{{\R^n}} (\phi (\bx)-\phi (0))^2 f_\infty(|\bx|)\ d\bx  = \int_{S^{n-1}}\int_0^\infty \left(\phi ( \rho,\Theta) -\phi (\mathbf 0)\right )^2   \rho^{n-1} f_\infty(\rho) \,d\rho\, d\Theta.
 \]
 Thanks to the gradient theorem, for the (smooth) function $\phi (\bx)$ we have
  \be\label{grad1}
\phi (\rho, \Theta) - \phi (0) = \int_0^1 \frac{\partial}{\partial \rho} \phi  (t \rho,\Theta) \, \rho \, dt.
\ee
Equality \fer{grad1} shows that the angles $\theta_i$, $i = 1,2,\dots, n-1$, appear as constant, so that the gradient theorem, for any given $\Theta$, takes a scalar form involving only a variation in the radius $\rho$. 
By Jensen's inequality, \fer{grad1} implies
\[
\left(\phi (\rho, \Theta) - \phi (0) \right )^2 \leq \int_0^1 \left|\frac{\partial}{\partial \rho} \phi  (t \rho,\Theta)\right |^2 \, \rho^2 \, dt
\]
and since $\nabla \phi (\bx)= 0$ if $|\bx| \le R$ we get
\begin{equations}\label{id4}
 Var[\phi (\bX)] \leq &  \int_{S^{n-1}} \int_0^\infty f_\infty(\rho) \rho^2 \left( \int_0^1 \left|\frac{\partial}{\partial \rho}  \phi  (t\rho,\Theta)\right|^2 \, dt\right) \rho^{n-1} \,d\rho\, d\Theta \\
 = & \int_{S^{n-1}} \int_R^\infty \left[\rho^n f_\infty(\rho)\right] \left(\rho \int_0^1 \left|\frac{\partial}{\partial \rho} \phi  (t \rho,\Theta)  \right|^2\, dt  \right) \,d\rho\, d\Theta.
\end{equations}
If we consider the $n$--dimensional equation \eqref{staz-2}
	\[
	\nabla[K(|{\bf x}|)f_\infty(|{\bf x}|)]+{\bf x} f_\infty(|{\bf x}|)=0,
	\]
	by multiplying by $\dfrac {{\bf x}}{|{\bf x}|}$ we get the scalar identity
	\be\label{scal}
	\frac{\partial }{\partial \rho} (K(\rho) f_\infty(\rho)) + \rho\, f_\infty(\rho) = 0.
	\ee
%
%
%
%
and so we obtain
\[
\rho^n \, f_\infty(\rho) = -  \frac{\partial }{\partial \rho} (K(\rho) \rho^{n-1} f_\infty(\rho)) + (n-1)  \frac{K(\rho)}{\rho^2} \rho^{n} f_\infty(\rho).
\]
 Hence, by resorting to inequality \fer{b1} we conclude that, for any $\rho \ge R$ it holds
 \be\label{k5}
 \frac 12 \rho^n \, f_\infty(\rho) \le -   \frac{\partial  }{\partial \rho}(K(\rho) \rho^{n-1} f_\infty(\rho)).
 \ee
 Making use of \fer{k5} into \fer{id4}, we get
    \begin{equations}\label{kkk}
  & \int_R^\infty \left[\rho^n f_\infty(\rho)\right] \left(\rho \int_0^1 \left|\frac{\partial}{\partial \rho} \phi  (t \rho,\Theta)  \right|^2\, dt\right) \,d\rho  \le \\ 
  &- 2  \int_R^\infty   \frac{\partial }{\partial \rho} (K(\rho) \rho^{n-1} f_\infty(\rho)) \left(\rho \int_0^1 \left|\frac{\partial}{\partial \rho} \phi  (t \rho,\Theta)  \right|^2\, dt\right) \,d\rho  \le   \\
  &2 \int_R^\infty K(\rho) \rho^{n -1} f_\infty(\rho)\frac{\partial}{\partial\rho} \left(\rho \int_0^1 \left|\frac{\partial}{\partial \rho} \phi  (t \rho,\Theta)  \right|^2\, dt\right) \,d\rho.
   \end{equations}
In \fer{kkk} we used  integration by parts to get the last line. Indeed here the border term satisfies
 \be\label{ss}
-  \left. \rho^{n-1}K(\rho)  f_\infty(\rho) \left(\rho \int_0^1 \left|\frac{\partial}{\partial \rho} \phi  (t \rho,\Theta)  \right|^2\, dt\right)    \right|_{R}^{+\infty} \le 0.
 \ee
Inequality \fer{ss} is a consequence of the fact that 
\[
\left.\int_0^1 \left|\frac{\partial}{\partial \rho} \phi  (t \rho,\Theta)  \right| ^2\, dt \right|_{\rho=R}=0.
\]
 The proof is completed by observing that, for any given function $\psi(\rho t)$ one has the identity
 \be\label{uno}
  \frac{\partial}{\partial \rho} \left[ \rho\, \psi(\rho  t)\right] =  \frac{\partial}{\partial t} \left( t\, \psi(\rho t )\right).
 \ee 
 Clearly, \fer{uno} implies
 \[
\begin{aligned}
& \frac{\partial}{\partial \rho} \left[\rho \int_0^1 \left|\frac{\partial}{\partial \rho} \phi  (t \rho,\Theta)  \right|^2\, dt \right] = \int_0^1 \frac{\partial}{\partial \rho} \left[\rho \left|\frac{\partial}{\partial \rho} \phi  (t \rho,\Theta)\right|^2 \right] \, dt = \\
&\int_0^1 \frac{\partial}{\partial t} \left[ t \left|\frac{\partial}{\partial \rho} \phi  (t \rho,\Theta)\right|^2 \right] \, dt  = 
\left|\frac{\partial \phi  (\rho,\Theta)}{\partial \rho}\right|^2.
\end{aligned}
 \]
Replacing into \fer{kkk}, and going back to cartesian coordinates gives the result for $\phi $
\[
\begin{split}
&Var[\phi (\bX)]  =  \int_{\R^n}  \left(\phi (\bx) - \int_{\R^n} \phi  (\by) f_\infty(|\by|) \, d\by\right)^2 f_\infty(|\bx|)\, d\bx\\
 &\leq   2  \int_{S^{n-1}}\int_R^\infty K(\rho) \rho^{n -1} f_\infty(\rho) \left|
 \frac{\partial  }{\partial \rho}\phi  (\rho,\Theta)\right|^2 \,d\rho \, d\Theta\\
 &\leq   2  \int_{S^{n-1}}\int_R^\infty K(\rho) \rho^{n -1} f_\infty(\rho) \left|\nabla \phi  (\rho,\Theta)\right|^2 \,d\rho \, d\Theta\\
 &=
   2  \int_{{\R^n}} K(|\bx|)  \left|\nabla \phi(\bx) \right|^2  f_\infty(|\bx| )\, d\bx.
\end{split}
\]
 \end{proof}
 
 \begin{rem} The previous proof can be easily adapted to the case in which $f_\infty(|\bx|)$ is a probability density supported in the ball $B_a = \{\bx: |\bx| \le a\}$, with $a$ being a positive constant. In this case, condition \fer{b1} is replaced by  the following: There exists a bounded value $0< R< a$ such that, for any $r$ such that $R\le r \le a$, the weight function $K(r)$, $r > 0$ satisfies  inequality \fer{b1}. An important example of this situation is provided by the class of Barenblatt densities defined in \fer{Bar-n}. 
 \end{rem}


\begin{rem} If we consider an isotropic probability density $f_\infty(|\bx|)$, the one-dimensional function $f_\infty(\rho)$ defined on $(0,+\infty)$ satisfies two differential equations. On one hand, if we consider the $n$--dimensional equation \eqref{staz-2}
	\[
	\nabla[K(|{\bf x}|)f_\infty(|{\bf x}|)]+{\bf x} f_\infty(|{\bf x}|)=0,
	\]
	we get the scalar identity
	\[
	\frac{\partial }{\partial \rho} (K(\rho) f_\infty(\rho)) + \rho\, f_\infty(\rho) = 0
	\]
	as we remarked in the previous proof.
	On the other hand, since  $\int_{{\R^n}} f_\infty(|{\bf x}|) d{\bf x}=1$, we get 
	\[
	\int_{0}^{+\infty} \sigma_n f_\infty(\rho) \rho^{n-1} \, d\rho =1.
	\]
	So, the function $f(\rho)= \sigma_n f_\infty(\rho) \rho^{n-1} $ is a probability density on $(0,+\infty)$ and can be considered as a steady state of many different Fokker--Planck equations as in \eqref{ste1}, namely
	\begin{equation}\label{star}
	\frac{\partial }{\partial \rho}
	\left(P(\rho)f(\rho)\right )+ 
	\,  Q(\rho) f(\rho) = 0.
	\end{equation}
	As a consequence, $f(\rho)$ satisfies many different weighted Poincaré inequalities as in \eqref{ch-gen} with a weight function $w$ associated to the diffusion and drift coefficients $P$ and $Q$ of the equation as in \eqref{mag}, namely
	\[
	w(\rho)= \dfrac{P(\rho)}{Q'(\rho)}, \quad \rho >0.
	\]
	In Examples \ref{sec:examples} in Section \ref{sec:poincare} we exploit the Poincaré inequalities coming from the Fokker--Planck equation \eqref{star}. On the contrary, in the proof of Proposition \ref{C}, we take advantage of the expression  \eqref{scal}.

\end{rem}


\section{Another hybrid inequality for isotropic densities}\label{sec:hybrid}
We are now presenting a refined  inequality for a $n$-dimensional isotropic density where we sharpen the contribution of the weight for large scales.  
Unfortunately, we are not able to get a genuine Poincaré--type inequality, due to an extra term on the right hand side. Nevertheless, this result applies without the restriction on the test function to be supported outside a ball.
\begin{thm}\label{H}
Let $\bX$ be a random vector distributed with density $f_\infty(|\bx|)$, with $\bx \in \R^n$ where the probability density function $f_\infty(|\bx|)$ satisfies the differential equation \eqref{staz-2}
\[
\nabla [K(|\bx|) f_\infty(|\bx|)] + \bx f_\infty(|\bx|) = 0.
 \]
 Suppose that there exists a real value $R>0$ such that, for any $r \ge R$, the weight function $K(r)$, $r \ge 0$ satisfies the inequality \eqref{b1}, namely
\[
(n-1) \frac{K(r)}{r^2} \le \frac 12 
\]
 and suppose moreover that  $f(\rho) =\sigma_n \rho^{n-1} f_\infty(\rho)$ satisfies the one-dimensional  Poincaré inequality \fer{Poi} on ${\it I}=(0,+\infty)$, with weight $w(\rho)$.
  Then, for any smooth function $\phi$  satisfying 
\be\label{smooth}
|\phi(\bx)| + |\nabla \phi(\bx)| \le C <\infty
\ee 
with $C=C(\phi) >0$, we have 
	\begin{equation}\label{chernoff-gen}
	\begin{gathered}
	\int_{\R^n} \phi^2(\bx)f_\infty(|\bx|) \, d\bx -  \left(\int_{\R^n} \phi(\by) f_\infty(|\by|) \, d \by\right )^2 \leq  \\
	C\left( \int_{\R^n} W(|\bx|) \left|\nabla\phi(\bx)\right |^2 f_\infty(|\bx|) \, d \bx + c(R)  \int_{|\bx|=R}  \left| \nabla \phi (\bx)\right |^2 f_\infty(|\bx|) d\sigma (\bx) \right )
	\end{gathered}
	\end{equation}
where 
\be\label{peso2}
W(|\bx|)= \max(w(|\bx|), |\bx|^2) \chi_{|x|\leq R} + K(|\bx|) \chi_{|x|>R},
\ee
$C>0$ and  $c(R)$ is a positive constant depending only on $R$.
\end{thm}
\begin{rem} The inequalities may then be classically extended to the class of all locally Lipschitz functions \cite{BL}. 
\end{rem}
\begin{rem} We underline that  in this theorem the weight function $W$ behaves at infinity as $K(|\bx|)$. If $K$ grows less than $|\bx|^2$ this result sharpens the previous Theorem \ref{Cher}, as long as the behavior at infinity of the weight is concerned. Nevertheless, we pay a prize in the extra term appearing in the right hand side.
\end{rem}

\noindent {\it Proof of Theorem \ref{H}.}
Let $\phi(\bx) = \phi(\rho, \Theta)$ where $(\rho, \Theta) = (\rho, \theta_1, \dots, \theta_{n-1})$ where $\rho =|\bx| \in(0,+\infty)$, $\theta_i\in (0,\pi)$, for any $i=1,\dots, n-2$ and  $\theta_{n-1}\in (0,2\pi)$. Moreover let $R>0$ as in the hypotheses of Theorem \ref{H} and let $0 <\e\ll R $.
Let $\bX$ be a random vector distributed with density $f_\infty(|\bx|)$, with $\bx \in  \R^n$. Then, if $\phi$ has finite variance, let us write in hyperspherical coordinates
\[
\begin{split}
& Var[\phi(\bX)]= \int_{\R^n} \left( \phi(\bx) - \int_{\R^n} \phi(\by) f_\infty(|\by|) \, d\by \right)^2 f_\infty(|\bx|)\, d\bx=\\
& \int_{S^{n-1}}\int_0^\infty \left(\phi(\rho, \Theta) - \int_{S^{n-1}}\int_0^\infty \phi(y, \eta) f_\infty(y) dy d \eta\right )^2 f_\infty (\rho) d\rho  d \Theta .
\end{split}
\]
We split now $\phi (\rho, \Theta)= \phi_\e^+(\rho, \Theta)+\phi_\e^-(\rho, \Theta)+\psi_\e(\rho, \Theta)$ where the three functions will be defined in a while.
We can estimate
\[
\begin{split}
& \int_{S^{n-1}}\int_0^\infty \left(\phi(\rho, \Theta) - \int_{S^{n-1}}\int_0^\infty \phi(y, \eta) f_\infty(y) dy d \eta\right )^2 f_\infty (\rho) d\rho  d \Theta =\\
 &\int_{S^{n-1}}\int_0^\infty \left[ 
 \left(\phi_\e^+(\rho, \Theta) - \int_{S^{n-1}}\int_0^\infty \phi_\e^+(y, \eta) f_\infty(y) dy d \eta\right )\right . \\
 &\qquad +  \left(\phi_\e^-(\rho, \Theta) - \int_{S^{n-1}}\int_0^\infty \phi_\e^-(y, \eta) f_\infty(y) dy d \eta\right )\\
 &\qquad +  \left . \left(\psi_\e(\rho, \Theta) - \int_{S^{n-1}}\int_0^\infty \psi_\e(y, \eta) f_\infty(y) dy d \eta\right )\right]^2 f_\infty (\rho) d\rho  d \Theta  \le \\
 &4 \left\{ 
 \int_{S^{n-1}}\int_0^\infty \left(\phi_\e^+(\rho, \Theta) - \int_{S^{n-1}}\int_0^\infty \phi_\e^+(y, \eta) f_\infty(y) dy d \eta\right )^2 f_\infty (\rho) d\rho  d \Theta  \right .\\
 &\qquad + \int_{S^{n-1}}\int_0^\infty \left(\phi_\e^-(\rho, \Theta) - \int_{S^{n-1}}\int_0^\infty \phi_\e^-(y, \eta) f_\infty(y) dy d \eta\right )^2 f_\infty (\rho) d\rho  d \Theta \\
 &\qquad + \left.\int_{S^{n-1}}\int_0^\infty 	\left(\psi_\e(\rho, \Theta) - \int_{S^{n-1}}\int_0^\infty \psi_\e(y, \eta) f_\infty(y) dy d \eta\right )^2 f_\infty (\rho) d\rho  d \Theta  
 \right \}\\
 & = 4 \left( Var[\phi_\e^+(\bX)] +Var[\phi_\e^-(\bX)]+ Var[\psi_\e(\bX)]\right ).
 \end{split}
 \]
Consequently, for any given $0<\e \ll R$ we can evaluate separately the variances of $\phi_\e^+$,  $\phi_\e^-$ and $\psi_\e$ with respect to $f_\infty$, gather the results and let $\e \to 0$ in order to get a global estimate.\\

We define now the three functions $\phi_\e^+$,  $\phi_\e^-$ and $\psi_\e$.
Let us consider the interval $(R-\e,R+\e)$ and let us set 
\begin{equation}\label{g+}
g_\e^+(\rho) = \begin{cases}
0  & \text{if } \rho \le R-\e, \\
P_\e^+(\rho)   & \text{if } R-\e < \rho < R, \\
1  & \text{if } \rho \ge  R.
\end{cases}
\end{equation}
In definition \fer{g+} the function $P_\e^+ (\rho)$ is the third order polynomial in $\rho$ connecting $C^1$ the points $(R- \e, 0)$  with  $(R,1)$, so that $g_\e^+(\rho)$ is $C^1((0,+\infty))$ with respect to $\rho$. Clearly, $P_\e^+ (\rho)$ is increasing with respect to $\rho$ in the interval $R-\e < \rho < R$.\\
Likewise, let us set
\begin{equation}\label{g-}
g_\e^-(\rho) = \begin{cases}
1 &\text{if } \rho \le  R, \\
P_\e^-(\rho)    &\text{if } R< \rho < R+\e, \\
0  &\text{if } \rho \ge R+\e.\\
\end{cases}
\end{equation}
In definition \fer{g-} the function $P_\e^- (\rho)$ is the third order polynomial connecting $C^1$ the points $(R,1)$ with $(R+\e, 0)$, so that $g_\e^-(\rho)$ is $C^1((0,+\infty))$ with respect to $\rho$. Clearly, $P_\e^- (\rho)$ is decreasing with respect to $\rho$ in the interval $R<\rho < R+\e$. We remark that, by construction, for any given $x \in (0,R)$ one has 
\begin{equation}\label{cond5}
\left\{
\begin{aligned}
& g_\e^+(R-x)= g_\e^-(R+x), \\
&\left. \frac{\partial}{\partial \rho}g_\e^+(\rho)\right|_{\rho = R-x} = - \left. \frac{\partial}{\partial \rho}g_\e^-(\rho)\right|_{\rho = R+x}. \\
\end{aligned}
\right.
\end{equation}
Let us finally define the approximations $\phi_{\e}^+$ and $\phi_{\e}^-$ of $\phi$ through the radial derivatives
\begin{equation}\label{phi}
\begin{aligned}
&\frac{\partial}{\partial\rho}\phi_{\e}^+(\rho,\Theta) =g_\e^+(\rho) \frac{\partial}{\partial\rho}\phi(\rho,\Theta)  \\
&\frac{\partial}{\partial\rho}\phi_{\e}^-(\rho,\Theta) = g_\e^-(\rho)\frac{\partial}{\partial\rho}\phi(\rho,\Theta).
\end{aligned}
\end{equation}
The functions  $\phi_{\e}^+$ and $\phi_{\e}^-$ are defined up to a function of $\Theta$.
Integrating by parts between $0$ and $\rho >0$, and setting  $\phi_\e^+(0, \Theta)= 0$ for all $\Theta$, we obtain
\[
\begin{split}
&\phi_\e^+(\rho, \Theta) = \int_0^\rho g_\e^+(y) \frac{\partial}{\partial y }\phi(y, \Theta)\, dy = \\
&\left[ g_\e^+(y) \phi(y, \Theta) \right]_0^\rho - \int_0^\rho \phi(y, \Theta)\frac{\partial}{\partial y }g_\e^+(y)\, dy.
\end{split}
\]
Hence, thanks to \fer{g+} we get
\begin{equation}\label{ff+}
\phi_\e^+(\rho, \Theta) =\begin{cases}
0 &\text{if } \rho \le R-\e, \\
P_\e^+(\rho) \phi(\rho,\Theta)  - {\displaystyle \int_{R-\e}^\rho }\phi(y, \Theta)\dfrac{\partial}{\partial y }P_\e^+(y)\, dy  & \text{if } R-\e < \rho < R, \\
 \phi(\rho,\Theta) -   {\displaystyle \int_{R-\e}^{R} }\phi(y, \Theta)\dfrac{\partial}{\partial y }P_\e^+(y)\, dy  & \text{if } \rho \ge  R.
\end{cases}
\end{equation}
Likewise we have
\begin{equations}\nonumber
&\phi_\e^-(\rho, \Theta) -  \phi_\e^-(0, \Theta) = \int_0^\rho g_\e^-(y) \frac{\partial}{\partial y }\phi(y, \Theta)\, dy = \\
&\left[ g_\e^-(y) \phi(y, \Theta) \right]_0^\rho - \int_0^\rho \phi(y, \Theta)\frac{\partial}{\partial y }g_\e^-(y)\, dy.
\end{equations}
Thus, setting $\phi_\e^-(0, \Theta)=\phi(0)$ for all $\Theta$,   \fer{g-} implies
\begin{equation}\label{ff-}
\phi_\e^-(\rho, \Theta) =
\begin{cases}
\phi(\rho,\Theta)  & \text{if } \rho \le R, \\
P_\e^-(\rho) \phi(\rho,\Theta)   -  {\displaystyle \int_{R}^\rho} \phi(y, \Theta)\dfrac{\partial}{\partial y }P_\e^-(y)\, dy  & \text{if } R< \rho < R+\e, \\
   -   {\displaystyle \int_R^{R+\e} }\phi(y, \Theta)\dfrac{\partial}{\partial y }P_\e^-(y)\, dy & \text{if } \rho \ge R+\e.
\end{cases}
\end{equation}
Setting $\psi_\e(\rho, \Theta)= \phi(\rho, \Theta) -  \phi_\e^+(\rho, \Theta)- \phi_\e^-( \rho, \Theta)$ we get
\begin{equation}\label{psi}
\psi_\e(\rho, \Theta) =
\begin{cases}
0 & \text{if } \rho \le R-\e,\\
\phi(\rho,\Theta)\left (1- P_\e^+(\rho)\right)  +  {\displaystyle \int_{R-\e}^\rho} \phi(y, \Theta)\dfrac{\partial}{\partial y }P_\e^+(y)\, dy 
   & \text{if } R-\e< \rho \le R, \\
\begin{split}
 -P_\e^-(\rho) \phi(\rho,\Theta)  &+      {\displaystyle \int_{R-\e}^R} \phi(y, \Theta)\dfrac{\partial}{\partial y }P_\e^+(y)\, dy \\
 &\qquad  +  {\displaystyle \int_{R}^\rho} \phi(y, \Theta)\dfrac{\partial}{\partial y }P_\e^-(y)\, dy 
 \end{split}
  & \text{if } R< \rho < R+\e, \\
     {\displaystyle \int_{R-\e}^{R+\e} }\phi(y, \Theta) 
     \left( \dfrac{\partial}{\partial y }g_\e^+(y) +  \dfrac{\partial}{\partial y }g_\e^-(y)\right )  \, dy & \text{if } \rho \ge R+\e.
\end{cases}
\end{equation}
Then, we have the identity 
\[
\phi(\rho, \Theta)= \phi_\e^+ (\rho,\Theta) + \phi_\e^- (\rho,\Theta) + \psi_\e(\rho, \Theta) , 
\]
where the function $\phi_\e^+ (\rho,\Theta) $ is equal to zero if $\rho \le R-\e$. 
Let us examine in details the size of the term $ \psi_\e(\rho, \Theta) $. Owing to conditions \fer{smooth} and \fer{cond5},  and since  $P^+_\e(\rho)$ is nondecreasing, and such that $0\le P^+_\e(\rho)\le 1$
\[
\left| \int_{R-\e}^\rho \phi(y, \Theta)\frac{\partial}{\partial y }P_\e^+(y)\, dy\right| \le\int_{R-\e}^\rho |\phi(y, \Theta)|\frac{\partial}{\partial y }P_\e^+(y)\, dy\le  C \int_{R-\e}^\rho \frac{\partial}{\partial y }P_\e^+(y)\, dy \le  C.
\]
The same upper bound holds for the  integrals appearing in the definition of $\psi_\e$ for $R<\rho<R+\e$. Therefore $\psi_\e(\rho, \Theta) =0 $ for $\rho \le R-\e$ and 
\begin{equations}\nonumber
|\psi_\e(\rho, \Theta)| &\le 3C \cdot \chi \left[ R-\e < \rho \le R+\e\right] \\
&+\chi\left[ \rho > R+\e\right]\left|\int_{R-\e}^{R+\e}  \phi(y, \Theta) \left[\frac{\partial}{\partial y }g_\e^+(y) + \frac{\partial}{\partial y }g_\e^-(y)\right] \, dy\right|.
\end{equations}
Moreover, thanks to conditions \fer{cond5}
\begin{equations}\nonumber
&\left| \int_{R-\e}^{R+\e}  \phi(y, \Theta) \left[\frac{\partial}{\partial y }g_\e^+(y) + \frac{\partial}{\partial y }g_\e^-(y)\right] \, dy \right|\\
& = \left| \int_{0}^\e \left[ \phi(R-x, \Theta) - \phi(R+x, \Theta)\right] \left. \frac{\partial}{\partial y }P_\e^+(y)\right|_{y= R-x}\, dx \right| \\
& \le \int_{0}^\e \left|\phi(R-x, \Theta) - \phi(R+x, \Theta)\right| \left. \frac{\partial}{\partial y }P_\e^+(y)\right|_{y= R-x}\, dx \\
&= \int_{0}^\e \frac{ \left|\phi(R-x, \Theta) - \phi(R+x, \Theta)\right| }{|2x|} |2x| \left. \frac{\partial}{\partial y }P_\e^+(y)\right|_{y= R-x}\, dx \\
& \le \sup |\nabla \phi(\bx)|\int_{0}^\e |2x| \left. \frac{\partial}{\partial y }P_\e^+(y)\right|_{y= R-x}\, dx \le 2C \e.
\end{equations}
Finally  $\psi_\e(\rho, \Theta) =0 $ for $\rho \le R-\e$ and 
\begin{equation}\label{bo1}
|\psi_\e(\rho, \Theta)| \le C\left( 3\chi \left[ R-\e < \rho< R+\e\right] +  2C\e 
\chi\left[ \rho > R+\e\right] \right).
\end{equation}
Inequality \fer{bo1} asserts that $|\psi_\e(\rho, \Theta)|$ is bounded by two constant functions, the first different from zero only on a set of measure $2\e$, the second bounded by $2C\e$ on the set $\rho > R+\e$. This clearly implies that 
\[
\lim_{\e\to 0} Var[\psi_\e( \bX)] =0.
\]
Let us come to $Var[\phi_\e^-(\bX)]$  and   $Var[\phi_\e^+(\bX)]$ and let us consider first $Var[\phi_\e^-(\bX)]$. Thanks to \eqref{Chern-sferiche}, we obtain
\begin{equation*}
Var[\phi_\e^-(\bX)] \le  c \int_0^\infty \int_{S^{n-1}} 	\left( w(\rho) \left| \frac{\partial}{\partial\rho}\phi_{\e}^-(\rho,\Theta) \right |^2 + \sum_{i=1}^{n-1}\left| \frac{\partial}{\partial \theta_i}\phi_{\e}^-(\rho,\Theta) \right |^2\right ) \rho^{n-1} f_\infty(\rho) d \Theta d\rho
\end{equation*}
where $w$ is the weight appearing in the one-dimensional Poincaré inequality relative to the density $\sigma_n \rho^{n-1} f_\infty(\rho)$ and $c>0$ is a suitable constant.
Let us split the integral in $d\rho$  into
\[
\int_0^\infty d\rho= \int_0^R d\rho+ \int_R^{R+\e}  d\rho+ \int _{R+\e}^\infty d\rho.
\]
Thanks to definition \eqref{ff-} we have
\[
\begin{split}
&\int_0^R \int_{S^{n-1}} \left( w(\rho) \left| \frac{\partial}{\partial\rho}\phi_{\e}^-(\rho,\Theta) \right |^2 + \sum_{i=1}^{n-1}\left| \frac{\partial}{\partial \theta_i}\phi_{\e}^-(\rho,\Theta) \right |^2\right ) \rho^{n-1} f_\infty(\rho) d \Theta d\rho\\
&=\int_0^R \int_{S^{n-1}} \left( w(\rho) \left| \frac{\partial}{\partial\rho}\phi(\rho,\Theta) \right |^2 + \sum_{i=1}^{n-1}\left| \frac{\partial}{\partial \theta_i}\phi(\rho,\Theta) \right |^2\right ) \rho^{n-1} f_\infty(\rho) d \Theta d\rho.
\end{split}
\]
Recalling that  $ \left|\dfrac{\partial \phi}{\partial\rho}\right|^2 \leq |\nabla \phi|^2 $ and $ \left|\dfrac{\partial \phi}{\partial\theta_i}\right|^2  \leq   \rho^2 |\nabla \phi|^2, \,  i=1, \dots, n-1$, we obtain
\begin{equation*}
\begin{split}
&\int_0^R \int_{S^{n-1}} \left( w(\rho) \left| \frac{\partial}{\partial\rho}\phi(\rho,\Theta) \right |^2 + \sum_{i=1}^{n-1}\left| \frac{\partial}{\partial \theta_i}\phi(\rho,\Theta) \right |^2\right ) \rho^{n-1} f_\infty(\rho) d \Theta d\rho \\
&\le n \int_0^R \int_{S^{n-1}} \max (w(\rho), \rho^2)  \left| \nabla \phi( \bx(\rho,\Theta)) \right |^2 \rho^{n-1} f_\infty(\rho) d \Theta d\rho\\
&= n \int_{|x| \leq R} \max (w(|\bx|), |\bx|^2)  \left| \nabla \phi( \bx) \right |^2  f_\infty(\bx) d\bx\\
&= n \int_{\R^n} \max (w(|\bx|), |\bx|^2) \chi_{|\bx| \le R}  \left| \nabla \phi( \bx) \right |^2  f_\infty(\bx) d\bx.
\end{split}
\end{equation*}
Now, let us estimate $\int_R^{R+\e} d \rho$. By definitions \eqref{phi} and  \eqref{ff-}, we have for $R<\rho<R+\e$ and for $i=1, \dots, n-1$
\begin{enumerate}[i)]
\item ${\displaystyle \frac{\partial}{\partial \rho}\phi_{\e}^-(\rho,\Theta)  =  P_\e^-(\rho)}\dfrac{\partial}{\partial \rho}\phi(\rho,\Theta)$
\item  ${\displaystyle \frac{\partial}{\partial \theta_i}\phi_{\e}^-(\rho,\Theta)  = P_\e^-(\rho)  \frac{\partial}{\partial \theta_i}\phi(\rho,\Theta) - \int_R^\rho 
\frac {\partial}{\partial \theta_i} \phi(y,\Theta) \frac {\partial}{\partial y} P_\e^-(y) dy}$.
\end{enumerate}
So, 
\begin{equation*}
\begin{split}
&\int_R^{R+\e} \int_{S^{n-1}} \left( w(\rho) \left| \frac{\partial}{\partial\rho}\phi_{\e}^-(\rho,\Theta) \right |^2 + \sum_{i=1}^{n-1}\left| \frac{\partial}{\partial \theta_i}\phi_{\e}^-(\rho,\Theta) \right |^2\right ) \rho^{n-1} f_\infty(\rho) d \Theta d\rho\\
&\le c  \int_R^{R+\e} \int_{S^{n-1}} \left( 
w(\rho) \left| \frac{\partial \phi}{\partial\rho}(\rho,\Theta) \right |^2 + \right . \\
&\qquad\qquad  \left.\sum_{i=1}^{n-1}\left[\left| \frac{\partial \phi}{\partial \theta_i}(\rho,\Theta) \right |^2 + 
\left( \int_R^\rho \left| \frac {\partial \phi}{\partial \theta_i}(y,\Theta)\right |\left(- \frac {\partial P_\e^-}{\partial y} (y)\right ) dy \right )^2 \right] 
\right ) \rho^{n-1} f_\infty(\rho) d \Theta d\rho\\
& \le 
c  \int_R^{R+\e} \int_{S^{n-1}} \left( 
w(\rho) \left| \frac{\partial \phi }{\partial\rho}(\rho,\Theta) \right |^2 +  
\sum_{i=1}^{n-1}\left[
\left| \frac{\partial \phi}{\partial \theta_i}(\rho,\Theta) \right |^2 + 
\sup_{R<\rho<R+\e} \left| \frac {\partial \phi} {\partial \theta_i}(\rho,\Theta)\right |^2 
\right] 
\right )\\
&\qquad\qquad\qquad \qquad \times  \rho^{n-1} f_\infty(\rho) d \Theta d\rho\\
& =A_\e.
\end{split}
\end{equation*}
Due to hypotesis \eqref{smooth} on the smoothness of $\phi$ and to the continuity of $w$, we can deduce by dominated convergence that 
\[
\lim_{\e \to 0} A_\e =0.
\]
Last, let us come to  $\int_{R+\e}^\infty d\rho$. For $\rho \ge R+\e$ the function $\phi_\e^-$ does not depend on $\rho$ so $\frac{\partial}{\partial \rho}\phi_\e^-(\rho, \Theta)=0$. Moreover,  for $\rho \ge R+\e$ and $i=1,\dots, n-1$ we have
\[
\frac{\partial}{\partial \theta_i} \phi_\e^-(\rho, \Theta) = - \int_R^{R+\e}  \frac{\partial}{\partial \theta_i} \phi(y, \Theta) \frac{\partial}{\partial y} P_\e^-(y) dy.
\]
So
\begin{equation*}
\begin{split}
&\int_{R+\e}^\infty \int_{S^{n-1}} 	\sum_{i=1}^{n-1}\left| \frac{\partial}{\partial \theta_i}\phi_{\e}^-(\rho,\Theta) \right |^2 \rho^{n-1} f_\infty(\rho) d \Theta d\rho\\
&= \int_{R+\e}^\infty \int_{S^{n-1}} 	\sum_{i=1}^{n-1}\left| 
-\int_R^{R+\e}  \frac{\partial}{\partial \theta_i} \phi(y, \Theta) \frac{\partial}{\partial y} P_\e^-(y) dy
 \right |^2 \rho^{n-1} f_\infty(\rho) d \Theta d\rho\\
 &\le  \int_{R+\e}^\infty \int_{S^{n-1}} 	\sum_{i=1}^{n-1}\left[
\int_R^{R+\e} \left| \frac{\partial}{\partial \theta_i} \phi(y, \Theta)\right | \left(-\frac{\partial}{\partial y} P_\e^-(y)\right ) dy
 \right ]^2 \rho^{n-1} f_\infty(\rho) d \Theta d\rho\\
 &\le   \int_{R+\e}^\infty \int_{S^{n-1}} 	
 \sum_{i=1}^{n-1}
 \sup_{R\le \rho \le R+\e} \left| \frac{\partial \phi}{\partial \theta_i} (\rho, \Theta)\right |^2 \rho^{n-1} f_\infty(\rho) d \Theta d\rho\\
 &= \int_{R+\e}^\infty
  \sum_{i=1}^{n-1}
 \left( \int_{S^{n-1}} 	
 \sup_{R\le \rho \le R+\e} \left| \frac{\partial \phi}{\partial \theta_i} (\rho, \Theta)\right |^2  d\Theta \right ) \rho^{n-1} f_\infty(\rho) d\rho \\
 &= 
  \left[\sum_{i=1}^{n-1}
 \left( \int_{S^{n-1}} 	
 \sup_{R\le \rho \le R+\e} \left| \frac{\partial \phi}{\partial \theta_i} (\rho, \Theta)\right |^2  d\Theta \right )\right ]
 \int_{R+\e}^\infty \rho^{n-1} f_\infty(\rho) d\rho \\
 &= B_\e
\end{split}
\end{equation*}
We get 
\[
\lim_{\e \to 0} \int_{S^{n-1}}  \sup_{R\le \rho \le R+\e} \left| \frac{\partial \phi}{\partial \theta_i} (\rho, \Theta)\right |^2  d\Theta = \int_{S^{n-1}}  \left| \frac{\partial \phi}{\partial \theta_i} (R, \Theta)\right |^2  d\Theta
\]
and so
\[
\begin{split} 
&\int_{S^{n-1}}  \left| \frac{\partial \phi}{\partial \theta_i} (R, \Theta)\right |^2  d\Theta \le \int_{S^{n-1}}  R^2 \left| \nabla \phi (\bx (R, \Theta))\right |^2  d\Theta\\
&= \frac {R^2}{R^{n-1}f_\infty(R)} \int_{S^{n-1}} \left| \nabla \phi (\bx (R, \Theta))\right |^2 R^{n-1} f_\infty(R) d\Theta \\
&=  \frac {R^3}{R^{n}f_\infty(R)} \int_{|\bx|=R}  \left| \nabla \phi (\bx)\right |^2  f_\infty(|\bx|) d\sigma (\bx)\\
&= c(R)  \int_{|\bx|=R}  \left| \nabla \phi (\bx)\right |^2  f_\infty(|\bx|) d\sigma (\bx).
\end{split}
\]
with $d\sigma(\bx)$ is the surface measure of the ball $|\bx|=R$ and $c(R)>0$ a constant depending only on $R$. 
By monotone convergence we get
\[
\begin{split}
\lim_{\e \to 0} B_\e &= c(R) \left( \int_{|\bx|=R}  \left| \nabla \phi (\bx)\right |^2 f_\infty(|\bx|) d\sigma (\bx) \right ) \left( \int_{R}^\infty \rho^{n-1} f_\infty(\rho) d\rho\right )\\
& \le 
c(R)  \int_{|\bx|=R}  \left| \nabla \phi (\bx)\right |^2 f_\infty(|\bx|) d\sigma (\bx).
\end{split}
\]
In the end
\[
Var[\phi_\e^-(\bX)]  \le   C \left(\int_{\R^n} \max (w(|\bx|), |\bx|^2)\chi_{|x| \leq R} \left| \nabla \phi( \bx) \right |^2 f_\infty(|\bx|) d\bx + A_\e + B_\e\right )\\
\]
and 
\begin{equation*}
\begin{split}
&\lim_{\e \to 0} Var[\phi_\e^-(\bX)] \\
& \le   C \left( \int_{\R^n} \max (w(|\bx|), |\bx|^2)\chi_{|x| \leq R} \left| \nabla \phi( \bx) \right |^2 f_\infty(|\bx|) d\bx + 
c(R)  \int_{|\bx|=R}  \left| \nabla \phi (\bx)\right |^2 f_\infty(|\bx|) d\sigma (\bx)\right )\\
& \le C \left(  \int_{\R^n} \max (w(|\bx|), |\bx|^2)\chi_{|x| \leq R}  \left| \nabla \phi( \bx) \right |^2 f_\infty(|\bx|) d\bx + 
c(R)  \int_{|\bx|=R}  \left| \nabla \phi (\bx)\right |^2 f_\infty(|\bx|) d\sigma (\bx)\right ).
\end{split}
\end{equation*}
We now estimate the last term $Var[\phi_\e^+(\bX)]$. Since $\phi_\e^+$ vanishes for $|\bx| \le  R-\e$, we can apply Proposition \ref{C} and obtain
\[
Var[\phi_\e^+(\bX)] \le 2  \int_{|\bx| \ge R-\e} K(|\bx|)  \left|\nabla \phi(\bx) \right|^2  f_\infty(|\bx| )\, d\bx
\]
where $K(|\bx|)$ is the weight related to the Fokker Planck equation \eqref{staz-2}.
Of course 
\[
\begin{aligned}
\lim_{\e \to 0} Var[\phi_\e^+(\bX)] & \leq   2  \int_{|\bx| \ge R} K(|\bx|)  \left|\nabla \phi(\bx) \right|^2  f_\infty(|\bx| )\, d\bx \\
&=   2  \int_{\R^n} K(|\bx|) \chi_{|\bx| \ge R} \left|\nabla \phi(\bx) \right|^2  f_\infty(|\bx| )\, d\bx.
\end{aligned}
\]
Since  
\[
Var[\phi(\bX)] \le 4 \lim_{\e \to 0} \left( Var[\phi_\e^+(\bX)] +Var[\phi_\e^-(\bX)]+ Var[\psi_\e(\bX)]\right )
\]
gathering all estimates we have obtained
\begin{equation}
	\begin{gathered}
	\int_{\R^n} \phi^2(\bx)f_\infty(|\bx|) \, d\bx -  \left(\int_{\R^n} \phi(\by) f_\infty(|\by|) \, d \by\right )^2 \leq  \\
	C \int_{\R^n} W(|\bx|) \left|\nabla\phi(\bx)\right |^2 f_\infty(|\bx|) \, d \bx + c(R)  \int_{|\bx|=R}  \left| \nabla \phi (\bx)\right |^2 f_\infty(|\bx|) d\sigma (\bx) 
	\end{gathered}
	\end{equation}
where $C>0$ is a suitable constant,
\be\label{peso2}
W(|\bx|)=  \max (w(|\bx|), |\bx|^2)\chi_{|\bx|\leq R} + K(|\bx|) \chi_{|\bx|>R}
\ee
and $c(R)$ is a positive constant depending only on $R$.
\hfill $\square$



\section{Appendix}
In order to prove Theorem \ref{giusto}, we show first the following.
\begin{prop}\label{buono} \cite{FPTT} Let $F(\bx,t)$ be the solution to equation \fer{quo} in $\MI \subseteq \R^n$. Then, if $\Psi(\bx)$ is a smooth function such that 
 \be\label{kk}
 \left| \Psi( \bx ) \right| \le c < \infty, \quad \bx \in \partial\MI
 \ee
for $t> 0$ the following equality holds
 \be\label{k3}
 \int_\MI f_\infty(\bx)\Psi(\bx) \frac{\partial F(\bx,t)}{\partial t} \, d\bx = -\int_\MI K(\bx)f_\infty(\bx) \nabla \Psi(\bx)\cdot \nabla F(\bx,t) \, d\bx.
 \ee
\end{prop}

\begin{proof}
The proof of \fer{k3} is immediate. Indeed, since $F(\bx, t)$ satisfies \fer{quo},  applying the divergence theorem in the volume $\MI$, we obtain 
 \begin{equations}\nonumber
 &\int_\MI f_\infty(\bx)\Psi(\bx) \frac{\partial F(\bx,t)}{\partial t}\, d\bx = \\
  & \int_\MI  \left[ \Psi(\bx) K(\bx)f_\infty(\bx)\Delta F(\bx,t)   - f_\infty (\bx) \Psi(\bx)(\mathbf {x-m})\cdot \nabla F(\bx,t) \right]\, d\bx = \\
&\int_{\{|\bx|=R\}}  \Psi(\bx) K(\bx)f_\infty(\bx) \nabla F(\bx,t)\cdot \mathbf n\, d\sigma (\bx)- \int_\MI \nabla F(\bx,t) \cdot \nabla\left( K(\bx)f_\infty(\bx) \Psi(\bx)\right)\, d\bx \, + \\
&- \int_\MI  f_\infty (\bx) \Psi(\bx)(\mathbf {x-m})\cdot \nabla F(\bx,t) \, d\bx= \\
& - \int_\MI K(\bx) f_\infty (\bx) \nabla F(\bx,t)\cdot\nabla\Psi(\bx) \, d\bx \, + \\
&- \int_\MI \Psi(\bx)\left[\nabla \left( K(\bx)f_\infty(\bx)\right) + f_\infty(\bx)(\mathbf {x-m}) \right]\cdot \nabla F(\bx,t)\, d\bx = \\
&- \int_IK(\bx) f_\infty  (\bx)  \nabla \Psi(\bx)\cdot \nabla F(\bx,t) \, d\bx.
 \end{equations}
Indeed, the border term vanishes in view of conditions \fer{kk} and  \fer{BCal2}, and we used \fer{stazio} on the last line to conclude.
\end{proof}

\noindent {\it Proof of Theorem \ref{giusto}:} Since the integral defining $\Theta(F(t))$ is uniformly bounded, we get
 \[
 \frac d{d t} \Theta (F(t)) =  \int_\MI f_\infty(\bx) \phi'(F(\bx,t)) \frac{\partial F(\bx,t)}{\partial t} \, d\bx.
 \]
Then, we apply Proposition \ref{buono} with $\Psi(\bx) =  \phi'(F(\bx,t))$ with fixed $t >0$.
\hfill $\square$

\subsection*{Acknowledgment}
This work has been written within the activities of GNFM (Gruppo Nazionale per la Fisica Matematica) and of GNAMPA (Gruppo Nazionale per l'Analisi Matematica, la Probabilit\`a  e le loro Applicazioni) of INdAM (Istituto Naziona\-le di Alta Matematica), Italy.
 The research was partially supported by the Italian  Ministry of University and Research (MUR) through the ``Dipartimenti di Eccellenza'' Programme 2023-2027 (CUP G43C22004580005) -- Dipartimento di Matematica dell'Università degli Studi di Milano and  through the MUR  project PRIN2022PNRR No.P2022Z7ZAJ. The authors states that there is no conflict of interest.

\end{document}